\documentclass{amsart}

\usepackage{amssymb,amsfonts, amsmath, amsthm}
\usepackage[all,arc]{xy}
\usepackage{enumerate}
\usepackage{mathrsfs}
\usepackage{mathtools}
\usepackage[mathscr]{euscript}
\usepackage{array}
\usepackage{appendix}

\usepackage{tikz}
\usepackage{tikz-cd}
\usepackage{mathpazo}
\usepackage{textcomp}
\usepackage{bbold}
\usepackage[bbgreekl]{mathbbol}

\usepackage[pagebackref, colorlinks=true, linkcolor=blue, citecolor = red]{hyperref}
\renewcommand*{\backref}[1]{}
\renewcommand*{\backrefalt}[4]{({%
		\ifcase #1 Not cited.%
		\or On p.~#2%
		\else On pp.~#2%
		\fi%
	})}

\usepackage[nameinlink,capitalise,noabbrev]{cleveref}
\crefname{subsection}{Subsection}{Subsection}

\usetikzlibrary{patterns}

\usepackage{geometry}
\usepackage{todonotes}

\DeclareMathAlphabet{\mathbbe}{U}{bbold}{m}{n}

\def\DDelta{{\mathbbe{\Delta}}}
\newcommand{\DD}{\DDelta}

\newcommand{\C}{\mathscr{C}}
\newcommand{\kC}{\mathfrak{C}}
\newcommand{\D}{\mathscr{D}}

\newcommand{\sL}{\mathscr{L}}

\newcommand{\N}{\mathscr{N}}

\newcommand{\kS}{\mathfrak{S}}

\newcommand{\comma}{,}

\newcommand{\set}{\mathscr{S}\mathrm{et}}
\newcommand{\cat}{\mathscr{C}\mathrm{at}}

\newcommand{\ds}{\displaystyle}
\newcommand{\sm}{\mathrm{sm}}

\newcommand{\Fun}{\mathrm{Fun}}
\newcommand{\Map}{\mathrm{Map}}
\newcommand{\Hom}{\mathrm{Hom}}

\newcommand{\id}{\mathrm{id}}

\newcommand{\Obj}{\mathrm{Obj}}

\newcommand{\QCat}{\mathscr{Q}\mathscr{C}\mathrm{at}}

\newcommand{\CSS}{\mathscr{C}\mathscr{S}\mathscr{S}}

\newcommand{\sint}{s\hspace{-0.04in}\int}
\newcommand{\sH}{s\mathscr{H}}

\newcommand{\ksS}{s\kS}

\newcommand{\kSeg}{\kS\mathrm{eg}}

\newcommand{\kCSS}{\kC\kS\kS}

\newcommand{\ksSet}{\mathrm{s}\mathfrak{S}\mathrm{et}}
\newcommand{\kssSet}{\mathrm{ss}\mathfrak{S}\mathrm{et}}
\newcommand{\ksssSet}{\mathrm{sss}\mathfrak{S}\mathrm{et}}

\newcommand{\ssL}{\mathrm{s}\sL}
\newcommand{\kT}{\mathfrak{T}}
\newcommand{\kI}{\mathfrak{I}}
\newcommand{\skT}{\mathrm{s}\mathfrak{T}}
\newcommand{\skI}{\mathrm{s}\mathfrak{I}}

\newcommand{\Fib}{\mathcal{F}\mathrm{ib}}
\newcommand{\LFib}{\mathcal{L}\Fib}

\newcommand{\LFibmin}{\LFib^{min}}

\newcommand{\Ree}{\mathcal{R}\mathrm{ee}}
\newcommand{\ReeLFib}{\Ree\LFib}

\newcommand{\Cart}{\mathcal{C}\mathrm{art}}
\newcommand{\coCart}{\mathrm{co}\Cart}

\newcommand{\SegcoCart}{\mathcal{S}\mathrm{egco}\mathcal{C}\mathrm{art}}

\newcommand{\Groth}{\mathscr{G}\mathrm{roth}}
\newcommand{\sSet}{\mathrm{s}\set}
\newcommand{\ssSet}{\mathrm{ss}\set}
\newcommand{\sssSet}{\mathrm{sss}\set}
\newcommand{\All}{\mathscr{A}\mathrm{ll}}
\newcommand{\sAll}{\mathrm{s}\All}

\newcommand{\Kan}{\mathscr{K}\mathrm{an}}

\newcommand{\kRee}{\mathfrak{R}\mathrm{ee}}
\newcommand{\Reemin}{\mathcal{R}\mathrm{ee}^{min}}

\newcommand{\kReemin}{\kRee^{min}}
\newcommand{\scMin}{\mathrm{s}\cMin}
\newcommand{\uKan}{\underline{\Kan}}
\newcommand{\uRee}{\underline{\Ree}}
\newcommand{\uKanmin}{\underline{\Kan}^{min}}
\newcommand{\uReemin}{\underline{\Ree}^{min}}

\newcommand{\cMin}{\mathcal{M}\mathrm{in}}
\newcommand{\kMin}{\mathfrak{M}\mathrm{in}}
\newcommand{\Lift}{\mathfrak{L}\mathrm{ift}}
\newcommand{\Sec}{\mathcal{S}\mathrm{ec}}
\newcommand{\sLift}{\mathrm{s}\Lift}

\newcommand{\uFun}{\underline{\Fun}}

\newcommand{\uGroth}{\underline{\Groth}}

\newtheorem{theone}[equation]{Theorem}

\newtheorem{lemone}[equation]{Lemma}
\newtheorem{propone}[equation]{Proposition}
\newtheorem{corone}[equation]{Corollary}

\theoremstyle{definition}
\newtheorem{defone}[equation]{Definition}
\newtheorem{exone}[equation]{Example}

\theoremstyle{remark}
\newtheorem{remone}[equation]{Remark}

\newtheorem{notone}[equation]{Notation}

\numberwithin{equation}{section}

\newtheoremstyle{TheoremNum}
{}{}              % space between body and thm
{\itshape}                      % Thm body font
{}                              % Indent amount (empty = no indent)
{\bfseries}                     % Thm head font
{.}                             % Punctuation after thm head
{ }                             % Space after thm head
{\thmname{#1}\thmnote{ \bfseries #3}}% Thm head spec
\theoremstyle{TheoremNum}
\newtheorem{thmn}{Theorem}

\makeatletter
\def\@seccntformat#1{%
  \expandafter\ifx\csname c@#1\endcsname\c@section\else
  \csname the#1\endcsname\quad
  \fi}
\makeatother

\title{A Model for the Higher Category of Higher Categories}

\author{Nima Rasekh}

\date{February 2022}

\begin{document}

\begin{abstract}
We use fibrations of complete Segal spaces as introduced in \cite{rasekh2017left,rasekh2021cartfibcss} to construct four complete Segal spaces: {\it Reedy fibrant simplicial spaces}, {\it Segal spaces}, {\it complete Segal spaces}, and {\it spaces}. Moreover, we show each one comes with a universal fibration that classifies {\it Reedy left fibrations}, {\it  Segal coCartesian fibrations}, {\it coCartesian fibrations} and {\it left fibrations} and prove these are representable fibrations in the sense of \cite{rasekh2017cartesian}. Finally, we use equivalences between quasi-categories and complete Segal spaces constructed in \cite{joyaltierney2007qcatvssegal,rasekh2021cartfibmarkedvscso} to present analogous constructions using fibrations of quasi-categories. 
\end{abstract}

\maketitle
\addtocontents{toc}{\protect\setcounter{tocdepth}{1}}

 \section{Introduction}\label{Introduction}
  
  \subsection{The (Higher) Category of (Higher) Categories}
  {\it Category theory} has been very effective in the study of a very diverse range of mathematical objects and their relation to each other. We can deduce various formal properties about different mathematical objects (such as the existence of free objects or preservation of universal properties) by using formal categorical results. A key illustration of these powerful methods is the study of sets via the category $\set$, which can be realized as the free cocompletion of the category with one object and so interesting properties, such as the fact that algebraic structures in sets are preserved by limits, follows formally \cite{maclane1998categories,riehl2017context}.
  
  This powerful perspective has been turned on category theory itself via the study of the (large) category of categories with objects small categories and morphisms functors, and we can similarly now deduce many valuable properties of categories, such as the construction of free categories, by analyzing properties of the category of categories.\footnote{We can in fact deduce all properties of categories by studying the $2$-category of categories, which is known as {\it formal category theory} \cite{gray1974formal}.} 
  
  While categories are a powerful tool in the study of classical mathematics, they are less suitable for objects that arise in homotopical mathematics. This starts with homotopy types of topological spaces or {\it Kan complexes} (which can be thought of as homotopical analogues of sets), but also {\it $A_\infty$-group}s \cite{stasheff1963ainftyspace} (which, up to homotopy, have a group structure) and even further {\it derived schemes} \cite{toen2014dag}. In order to effectively study such homotopical objects various notions of weak or homotopical categories have been developed, now known as {\it higher categories} or {\it $(\infty,1)$-categories} or simply {\it $\infty$-categories} \cite{bergner2010survey}. The most popular model are {\it quasi-categories} \cite{boardmanvogt1973qcats}, and other important models are {\it Kan enriched categories} \cite{dwyerkan1980simplocalization,bergner2007bergnermodelcat} and {\it complete Segal spaces} \cite{rezk2001css}, which are related to each other via various equivalences \cite{joyaltierney2007qcatvssegal,bergner2007threemodels,lurie2009htt}. These various models of $\infty$-categories give us the appropriate framework to study concepts such as homotopical algebra and derived geometry.
  
  The analogy with classical category theory would suggest that similar to the category of sets, there is an easily constructed higher category of spaces that can be studied using higher categorical methods. While there is such a higher category, the construction is by no means immediate. The situation gets worse when trying to construct the higher category of higher categories. It is in fact a cruel joke of higher categorical mathematics that the construction of the higher category of higher categories requires us to change models, making the construction quite complicated. 
  
  \subsection{Strict Categories and Nerves} \label{subsec:nerves}
  The easiest way to construct an $(\infty,1)$-category is via Kan enriched categories and so our first way to approach this problem is by constructing a Kan enriched category. Constructing the Kan enriched category of Kan complexes is fairly straightforward and has been known at least since work of Quillen \cite{quillen1967modelcats}. We can use a similar line of thinking to construct an $(\infty,1)$-category of $(\infty,1)$-categories. Indeed, following work of Rezk, the category of complete Segal spaces is in fact enriched over Kan complexes \cite{rezk2001css}. On the other hand the category of quasi-categories is not directly enriched over Kan complexes but rather over quasi-categories themselves, however, we can easily construct a Kan enriched category by taking the underlying Kan complexes of the mapping quasi-categories \cite{riehlverity2022elements}. On the other side, as suggested by the slogan above, there is no immediate way to enrich the category of Kan enriched categories over Kan complexes.
  
  While the construction of the Kan enriched category of spaces and $(\infty,1)$-categories might initially appear to be a satisfactory answer, we are quickly confronted with various challenges. First of all most of $(\infty,1)$-category theory has not been developed in the context of Kan enriched categories\footnote{In fact even the construction of a simple over-category can be a challenge in this setting \cite[Appendix A]{horev2022moorecat}}. Historically speaking most of higher category theory has been developed using quasi-categories \cite{joyal2008notes,joyal2008theory,lurie2009htt}. Moreover, in an effort to move beyond one specific model Riehl and Verity developed a new method to approach higher category theory model-independently via the notion of an {\it $\infty$-cosmos}, which can then in particular be applied to quasi-categories, but also complete Segal spaces, Segal categories and even {\it $1$-complicial sets}, but notably not Kan enriched categories \cite{riehlverity2022elements}. 
  
  This motivates a construction of the higher category of spaces and $(\infty,1)$-categories using more established models, such as quasi-categories and complete Segal spaces. One first approach might be to simply translate the construction from Kan enriched categories to these models using various {\it nerve} constructions, which take a Kan enriched categories to (bi)simplicial sets. In particular, we can use the {\it simplicial nerve} \cite{lurie2009htt} to construct quasi-categories out of Kan enriched categories. Similarly we can use the {\it Rezk nerve}\footnote{Called {\it classifying diagram} in \cite{rezk2001css}.} or its variant due to Barwick and Kan to construct complete Segal spaces out of Kan enriched categories \cite{barwickkan2012relativecategory,meier2016fibcat}. These constructions are theoretically very satisfying, however, are computationally very challenging. Indeed, the construction due to Rezk requires a fibrant replacement in the Reedy model structure \cite[Section 8]{rezk2001css}, which, while preserving the level-wise homotopy type of the complete Segal space, complete changes the point-set structure. The simplicial nerve \cite[Proposition 1.1.5.10]{lurie2009htt} and the nerve by Barwick-Kan \cite[Theorem 6.1]{barwickkan2012relativecategory} do not require such additional steps, however, are by definition far more complicated constructions.\footnote{Indeed, there are several papers dedicated to understanding the left adjoint $\mathfrak{C}$ to the simplicial nerve \cite{duggerspivak2011rigidification,riehl2011simpcatofqcat}.}
  
\subsection{Higher Categories via Fibrations} \label{subsec:higher cat fib}
 What makes the nerve constructions so complicated is the fact that maps of spaces (and functors of $\infty$-categories) are by definition strict functors and so we either need to use a very complicated nerve construction (such as the Rezk nerve or simplicial nerve) or use a naive construction and then apply a Reedy fibrant replacement, both with the goal of ``destrictifying" the functors in order to allow for the possibility of higher categorical pseudo-functors, whose functoriality only holds up to higher equivalences. Ideally, we could have directly constructed a Kan-enriched category of spaces or $\infty$-categories where the morphisms directly correspond to some notion of pseudo-functors. However, directly defining pseudo-functors would require specifying an infinite tower of data and so we need to find a way to circumvent this dilemma. We need to choose a notion of functor of spaces (and $\infty$-categories) that is by definition weaker, yet still manageable. Fortunately, there is already an excellent solution in the category theory literature: {\it Fibrations}.
 
 The idea of using fibrations as a replacement for functors goes back to work of Grothendieck and Bourbaki, who used the (now called) {\it Grothendieck fibrations} to study set-valued and category-valued functors \cite{grothendieck2003etalegroup}. In particular, Grothendieck opfibrations over a category $\C$ correspond to a pseudo-functor $\C \to \cat$. Following the definition of a pseudo-functor  \cite{benabou1967bicat}, a pseudo-functor $[0] \to \cat$, is the data of a category $\C$ and an automorphism that is naturally isomorphic to the identity. Similarly, a pseudo-functor $[1] \to \cat$ is the data of a functor $\C \to \D$, and choices of automorphisms of $\C$ and $\D$, that all interact with each other in the appropriate manner. These examples already illustrate that by taking Grothendieck opfibrations over the categories $[n]$, we obtain a much less rigid object than we would if we use the classical nerve, which is defined as $N\cat_n = \Fun([n],\cat)$ i.e. strict functors out of $[n]$.
 
 This philosophy expands to the $\infty$-categorical setting. We hence want to construct a quasi-category and complete Segal space of spaces and $\infty$-categories by choosing an appropriate notion of fibration over the representable diagrams (i.e. the appropriate analogue to $[n]$). This requires use to use the vast literature on fibrations of $\infty$-categories. Concretely, for many theories of $\infty$-categories (such as quasi-categories, complete Segal spaces and in fact every other $\infty$-cosmos) $\infty$-categorical functors from an $\infty$-category $\C$ into spaces correspond to {\it left fibrations}\footnote{Also called {\it groupoidal coCartesian fibration} over $\C$ in \cite{riehlverity2022elements}.}  \cite{debrito2018leftfibration,heutsmoerdijk2015leftfibrationi,rasekh2017left,cisinski2019highercategories}. Similarly, functors valued in $\infty$-categories are classified by {\it coCartesian fibrations} over $\C$ \cite{rasekh2021cartfibcss,rasekh2021cartfibmarkedvscso,ayalafrancis2020fibrations,lurie2009htt,riehlverity2022elements}. 
 
\subsection{Constructing Complete Segal Spaces via Fibrations}
 The goal of this paper is to make the intuition outlined in the previous section into precise mathematical statements. In particular, we prove that the bisimplicial set $\kS$, which is characterized level-wise as left fibrations of bisimplicial sets over the representables $\Delta[n,l]$ (\cref{not:delta dash dash}), satisfies the following conditions.

 \begin{thmn}[\ref{the:main spaces}]
 	There is an equivalence between the complete Segal space $\kS$ and the strict Segal category $N_\Delta\uKan$ (\ref{eq:nerve}). Moreover, we have a natural bijection $\LFib \cong \Hom(-,\kS)$ (\ref{eq:lfib}).
 \end{thmn}
 
 This corresponds to a similar result by Kazhdan and Varshavsky \cite{kazhdanvarshvsky2014yoneda} (\cref{rem:kazhdan}) and generalizes a result by Kapulkin and Lumsdaine \cite{kapulkinlumsdaine2021kanunivalent} (\cref{rem:kapulkin}). 

Having used left fibrations of simplicial spaces to construct the complete Segal space of spaces we next generalize our result in order to construct an $\infty$-category of $\infty$-categories. We in fact obtain a far more general result. Using the observation that every complete Segal space is a Reedy fibrant simplicial space, we first use the theory of {\it Reedy left fibrations} to construct the $\infty$-category of Reedy fibrant simplicial spaces.

\begin{thmn}[\ref{the:main simp spaces}]
	There is an equivalence between the complete Segal space $\ksS$ and the strict Segal category $N_\Delta\uRee$ (\ref{eq:map reedy}). Moreover, we have a bijection $\ReeLFib \cong \Hom(-,\ksS)$.
\end{thmn}

While the construction of the complete Segal space of simplicial spaces could have relevance in the study of certain type theories that arise in the work of Riehl and Shulman \cite{riehlshulman2017rezktypes} (as further discussed in \cref{rem:shulman}) our main focus here is to restrict to the sub-complete Segal space of (complete) Segal spaces.

\begin{thmn}[\ref{the:css}]
	We have a diagram of fully faithful functors of complete Segal spaces 
	$$\kCSS \hookrightarrow \kSeg \hookrightarrow \ksS,$$
	with $\kSeg$ having elements Segal spaces and $\kCSS$ complete Segal spaces. Moreover, we have bijections 
	$$\SegcoCart \cong \Hom_{\ssSet}(-,\kSeg),$$
	$$\coCart \cong \Hom_{\ssSet}(-,\kCSS).$$
\end{thmn}

 The existence of the desired complete Segal spaces of (simplicial) spaces and (complete) Segal spaces with the universal property outlined above implies the existence of universal fibrations, which is the focus of \cref{sec:universal}. We in particular establish that the {\it universal left fibration} is represented by the terminal object (\cref{the:lfib rep}) and that the {\it universal Reedy left fibration}, Segal coCartesian fibration and coCartesian fibration are represented, in the sense of \cite{rasekh2017cartesian}, by the cosimplicial object $\DD \to \cat_\infty$ taking $[n]$ to $n$ composable morphisms (\cref{the:ree lfib rep}/\cref{cor:universal fib css}), which has also been discussed in \cite{rasekh2017cartesian,stenzel2020comprehension} (\cref{rem:rep cart fib}).
 
 While most of our work focuses on constructing complete Segal spaces, in the last section we use equivalences constructed by Joyal and Tierney \cite{joyaltierney2007qcatvssegal} and its fibrational analogue \cite{rasekh2017left,rasekh2021cartfibmarkedvscso} to construct various quasi-categories, beginning with the quasi-category of spaces:
 
 \begin{thmn}[\ref{the:qcat v segal}]
 	The maps $\kT:\kS_{\QCat} \to i_1^*\kS$ (\ref{eq:kt}) and $\kI: i_1^*\kS \to \kS_{\QCat}$ (\ref{eq:ki}) are inverses of quasi-categories.
 \end{thmn} 

This hence recovers a result by Cisinski \cite{cisinski2019highercategories} (\cref{rem:cisinski}), who uses a similar idea to construct the quasi-category of spaces using left fibrations of simplicial sets. As a next step, we can further generalize the result to quasi-category of simplicial spaces.
 
 \begin{thmn}[\ref{the:qcat v segal s}]
 	The maps $\skT:\ksS_{\QCat} \to i_1^*\ksS$ and $\skI: i_1^*\ksS \to \ksS_{\QCat}$  (\ref{eq:skt}) are inverses of quasi-categories.
 \end{thmn}

 We can restrict this equivalence easily to (complete) Segal spaces to obtain further equivalences (\cref{cor:big diagram}).
  
\subsection{Background and Notation} \label{subsec:notation}
We will assume familiarity with standard category theory as can be found in \cite{maclane1998categories,riehl2017context}. Also some familiarity with complete Segal spaces \cite{rezk2001css} and quasi-categories \cite{rezk2017qcats} would be helpful. Moreover, we make extensive use of left and coCartesian fibrations as studied in \cite{rasekh2017left,rasekh2017cartesian,rasekh2021cartfibcss,rasekh2021cartfibmarkedvscso}, however, key results have been reviewed when necessary.

Throughout, $\DD$ denotes the simplex category. We both use the functor category and the set of functors and hence denote the functor category by $\uFun$, whereas the set of functors is $\Fun$. Moreover for a given functor $F: \C \to \D$, the {\it set} of objects in the category $\uFun(\C,\D)_{/F}$ is also denoted by $\Fun(\C,\D)_{/F}$.

Finally, we will have a category of small sets, denoted $\set^{\sm}$, and for every category $\C$ and functor $X:\C \to \set$, we use $(\Fun(\C,\set)_{/X})^{\sm}$ to denote the full subcategory of $\Fun(\C,\set)_{/X}$ with objects maps $\alpha:Y \to X$ with small fiber.

\subsection{Acknowledgments} 
 I would like to thank William Balderrama for many helpful conversations and in particular for making me aware of 
 \cite[Definition 5.2.3]{cisinski2019highercategories}. I also would like to thank Emily Riehl for many helpful conversations
 and in particular for making me aware of the non-functoriality of pullbacks. I also would like to thank Hoang Kim Nguyen for discussion regarding minimal fibrations. Moreover, I would like to thank Denis-Charles Cisinski for several helpful explanations regarding his results in \cite{cisinski2019highercategories}. Finally, I would also like to thank the referee for many helpful suggestions.
 
\section{Background $\&$ Technicalities} \label{sec:tech}
As explained in \cref{subsec:higher cat fib}, we want to construct a complete Segal space of spaces $\kS$ (which we will do in \cref{sec:spaces}) which is level-wise given by a set of left fibrations over representable objects. However, there are two significant theoretical challenges that we need to overcome. First of all, if we naively define $\kS_{nl}$ as the set of left fibrations over $\Delta[n,l]$ (\cref{not:delta dash dash}), then the functoriality needs to follow from pulling back left fibrations. However, a pullback is only determined up to isomorphisms, hence a functor $\DD^{op} \times \DD^{op} \to \set$ that takes each pair $([n],[l])$ to the set of left fibrations over $\Delta[n,l]$ and each morphisms to the pullback would only be {\it pseudo-functorial}. In order to avoid this problem we associate functors to our fibrations that can then be strictly composed, which is the goal of \cref{subsec:fib functor}. 

Having taking care of the pseudo-functoriality, we can in fact directly define a bisimplicial set $\kS$ with $\kS_{nl}$ given by left fibrations over $\Delta[n,l]$ and we would like to prove that this is in fact a complete Segal space. Here the next problem arises. In order to prove that $\kS$ is a complete Segal space we need to show that for every trivial cofibration in the complete Segal space model structure $i:A \hookrightarrow B$ the map 
$$i^*:\Hom(B,\kS) \to \Hom(A,\kS)$$
is surjective (\cref{rem:css lift}). As we will establish in \ref{eq:lfib} this is equivalent to $i^*:\LFib(B) \to \LFib(A)$ being surjective, which means we need to prove that every left fibration over $A$ can be obtained as a pullback of a left fibration over $B$. By \ref{eq:derived unit} it is immediate that every left fibration is a {\it homotopy} pullback of a left fibration over $B$, however we need a strict pullback. In order to guarantee we can obtain this strict pullback we need to review the theory of minimal fibrations, which is the goal of \cref{subsec:min fib}.

\begin{remone}
	There are alternative ways to the ones introduced in \cref{subsec:fib functor} to avoid the pseudo-functoriality of the pullback of fibrations. For example we can choose well-orderings on all the fibers of the fibrations as has been done in \cite[Subsection 2.1]{kapulkinlumsdaine2021kanunivalent} or one can make a choice of a collection of pullbacks as has been done in \cite[Subsection 2.1]{cisinski2019highercategories}. On the other hand all these sources also rely on minimal fibrations (similar to \cref{subsec:min fib}) to guarantee the existence of strict lifts. 
\end{remone}

\subsection{Simplicial Objects and Fibrations} \label{subsec:simplicial}
In this short section we give a quick review of the necessary simplicial objects, relevant notation and their fibrations. 

\begin{notone} \label{not:delta dash dash}
 Denote the category of simplicial sets by $\sSet$ with generators $\Delta[-]$. Moreover, use $\partial \Delta[-]$ to denote its boundary and $\Lambda[-]$ its horn \cite[Subsection I.3]{goerssjardine2009simplicialhomotopytheory}. Similarly, denote the category of bisimplicial sets by $\ssSet$ and the generators by $\Delta[-,-]$. Notice both categories are simplicially enriched and we denote the enrichment by $\Map(-,-)$ \cite[Subsection I.V]{goerssjardine2009simplicialhomotopytheory}.
\end{notone} 
 For a given bisimplicial set $X_{\bullet\bullet}$ we use the notation $X_n$ to denote the simplicial set $(X_n)_l = X_{nl}$. Our first important fibration of bisimplicial sets are {\it Reedy fibrations}.

\begin{defone} \label{def:reedy}
	A {\it Reedy fibration} is a map of bisimplicial set $p:Y \to X$ that has the right lifting property with respect to the maps 
	$$\Delta[n,0] \times \Lambda[0,l]_i \coprod_{\partial \Delta[n,0] \times \Lambda[0,l]_i} \partial \Delta[n,0] \times \Delta[0,l] \to \Delta[n,0] \times \Delta[0,l] \cong \Delta[n,l],$$
	for all $n, l \geq 0$ and $0 \leq i \leq l$. A Reedy fibration is moreover  {\it trivial} if it satisfies the right lifting property with respect to all inclusions of bisimplicial sets.
\end{defone}

The Reedy fibrancy condition is equivalent to the map $Y_n \to X_n \times_{M_nX} M_nY$ being a Kan fibration of simplicial sets, where $M_nY,M_nX$ are the {\it matching spaces}, which in the particular case of bisimplicial sets are given by the simplicial set $M_nX = \Map_{\ssSet}(\partial \Delta[n,0],X)$. Reedy fibrations are part of a model structure with cofibrations given by inclusions of bisimplicial sets and equivalences given by level-wise Kan equivalences. In particular all trivial Reedy fibrations as Reedy weak equivalences. See \cite[Subsection 2.4]{rezk2001css} and \cite[Subsection 5.1]{hovey1999modelcategories} for more details. Reedy fibrations can be used to define a prominent model of $(\infty,1)$-categories, complete Segal spaces, which are defined by Rezk \cite{rezk2001css} and proven to be models of $(\infty,1)$-categories in \cite{bergner2007threemodels,joyaltierney2007qcatvssegal,toen2005unicity}. 

\begin{defone} \label{def:css}
	A {\it complete Segal space} $W$ is a Reedy fibrant simplicial space that satisfies the following two conditions:
	\begin{itemize}
		\item {\bf Segal Condition:} The restriction map 
		$$W_n \to W_1 \times_{W_0} ... \times_{W_0} W_1$$
		is a Kan equivalence for all $n \geq 2$.
		\item {\bf Completeness Conditions:} The map $W_0 \to W_0 \times W_0 \times_{W_1 \times W_1} W_3$ is a Kan equivalence.
	\end{itemize}
\end{defone}

Complete Segal spaces are can be used to do $(\infty,1)$-category theory. Here we think of the objects as the set $W_{00}$ and for two objects $x,y$ in $W$, the {\it mapping space} is defined as 
\begin{equation} \label{eq:map}
	\Map_W(x,y) = \Delta[0] \times_{W_0 \times W_0} W_1.
\end{equation}
For more details regarding the category theory of complete Segal space see \cite[Section 5]{rezk2001css}. Complete Segal spaces are in fact fibrant objects in a model structure, the {\it complete Segal space} model structure on $\ssSet$, defined originally by Rezk \cite[Theorem 7.2]{rezk2001css}. 

\begin{remone}\label{rem:css lift}
 This in particular implies that a bisimplicial set $W$ is a complete Segal space if and only if for every trivial cofibration $i:A \to B$ in the complete Segal space model structure, the map $i^*:\Hom(B,W) \to \Hom(A,W)$ is surjective, meaning every map $f: A \to W$ lifts along $i$.
\end{remone}

 Next we have left fibrations.

\begin{defone} \label{def:left fib}
	A left fibration is a Reedy fibration of bisimplicial sets $p:L \to X$ that satisfies the right lifting property with respect to maps 
	$$\Delta[n,0] \times \partial \Delta[0,l] \coprod_{\Delta[0,0] \times \partial \Delta[0,l]} \Delta[0,0] \times \Delta[0,l] \to \Delta[n,0] \times \Delta[0,l] \cong \Delta[n,l],$$
	induced by the map $\{0\}:\Delta[0,0] \to \Delta[n,0]$ i.e. the morphism that corresponds to $0 \in \Hom(\Delta[0,0],\Delta[n,0]) = \Delta[n,0]_{00} = \{0,...,n\}$.
	By \cite[Proposition 3.7]{rasekh2017left} this is equivalent to being a Reedy fibration and for all $n$, the map 
	\begin{equation} \label{eq:left}
		L_n \to X_n \times_{X_0} L_0
	\end{equation} 
	being a trivial Kan fibrations. 
\end{defone} 

\begin{remone} \label{rem:lfib constant}
	If $X = \Delta[0,0]$ then the condition \ref{eq:left} implies $L$ is {\it homotopically constant}, meaning $L_n \simeq L_0$ for all $n \geq 0$, and so $L$ is (homotopically) uniquely determined by the Kan complex $L_0$.
\end{remone}

\begin{remone} \label{rem:strict left}
	Notice if $L \to X$ is a map of bisimplicial sets, such that $L_n \cong X_n \times_{X_0} L_0$, meaning $L_n$ is a strict pullback, then the Reedy fibrant replacement $\hat{L} \to X$ of $L$ over $X$ is both a Reedy fibration and $L_n \simeq X_n \times_{X_0} L_0$, meaning $\hat{L} \to X$ is a left fibration.
\end{remone}

Left fibrations and complete Segal space equivalences interact well with each other. Concretely, for every complete Segal space equivalence $i:A \to B$ and left fibration $L \to A$, there exists a left fibration $\hat{L} \to B$ and a homotopy pullback square
\begin{equation} \label{eq:derived unit}
	\begin{tikzcd}
		L \arrow[r] \arrow[d, twoheadrightarrow] \arrow[dr, phantom, "\ulcorner", very near start] & \hat{L} \arrow[d, twoheadrightarrow] \\
		A \arrow[r, "i"] & B
	\end{tikzcd}
\end{equation}
given by the derived unit of the Quillen equivalence constructed in \cite[Theorem 5.1]{rasekh2017left}, meaning $A$ is Reedy equivalent to $i^*B$.

\begin{remone}\label{rem:rfib}
	Analogous to the definition of left fibrations (\ref{eq:left}) we can also define right fibrations as Reedy fibrations $R \to X$ such that $R_n \to X_n \times_{X_0} R_0$ is an equivalence, this time induced by the map $\{n\}:\Delta[0,0] \to \Delta[n,0]$ \cite[Remark 4.24]{rasekh2017left}. 
	Right fibrations are completely determined by left fibrations. Indeed, let $(-)^{op}:\ssSet \to \ssSet$ be the automorphism induced by the unique non-trivial automorphism $\sigma \times \sigma$ from $\DD\times\DD$ to itself. Then a map $R \to X$ is a right fibration if and only if $R^{op} \to X^{op}$ is a left fibration.
\end{remone}

\subsection{Fibrations vs. Functors} \label{subsec:fib functor}
In this subsection we construct a precise way to translate between functors and fibrations to avoid the pseudo-functoriality that arises when using pullback (as discussed in the beginning of \cref{sec:tech}). We will start by reviewing basic facts regarding functors and fibrations as discussed in \cite{maclanemoerdijk1994topos} or \cite{johnstone2002elephanti,johnstone2002elephantsii}. 

Recall that a {\it discrete Grothendieck fibration} is a functor $p: \D \to \C$ such that for every morphism $f: c \to c'$ and chosen lift $d'$ in $\D$ (meaning $p(d')=c'$) there exists a {\it unique} $\hat{f}:\hat{c} \to \hat{c}'$ in $\D$ such that $p(\hat{f}) =f$. Following the convention from \cref{subsec:notation}, we use $\uGroth(\C)$ to denote the full subcategory of $\cat_{/\C}$ with objects discrete Grothendieck fibrations with small fibers and $\Groth(\C)$ for the large set of discrete Grothendieck fibrations with small fibers. It is well established that $\uGroth(\C)$ is equivalent to the functor category $\uFun(\C^{op},\set)$\footnote{In fact we have far more general results for Grothendieck fibrations \cite[Theorem B1.3.6]{johnstone2002elephanti}.}, however, we need an isomorphism and hence state the desired result explicitly.

\begin{lemone} \label{lemma:groth iso}
	Let $\C$ be a small category. There is an isomorphism of categories 
	\begin{center}
		\begin{tikzcd}
			\uFun(\C^{op},\set) \arrow[rr, "\int_\C", "\cong"', shift left=1.8]& & \uGroth(\C) \arrow[ll, "\Fib_\C", shift left=1.8]
		\end{tikzcd}
	\end{center} 
\end{lemone}

\begin{proof}
	We take the proof of the equivalence as given in \cite[Theorem 2.1.2]{loregianriehl2020fib} and show it is in fact an isomorphism. First we review the construction of the relevant functors. For a given functor $F: \C^{op} \to \set$, define $\int_\C F$ as the category with objects the disjoint union $\coprod_{c \in \Obj\C} F(c)$ and morphisms $f: x \in F(c) \to y \in F(d)$ if there exists a morphism $f: c \to d$ in $\C$ such that $F(f)(y) = x$. On the other hand, let $\Fib(P:\D \to \C): \C^{op} \to \set$ be the functor that takes $c$ to the fiber $\Fib(c) = P^{-1}(c)$, which is indeed a functor using the unique extension property of discrete Grothendieck fibrations. 
	
	Now, $\Fib\int_\C F$ is a functor with value the fiber of $\int_\C F$ over $c$, which by construction of the disjoint union is $F(c)$. On the other side, $\int_\C\Fib P$ is a functor over $\C$ with objects over an object $c$ in $\C$ given by the value $\Fib P(c)$, which by definition is just the fiber of $P$ over $c$. 
\end{proof}

This lemma has two important corollaries. Let $\Groth = \coprod_{\C \in \cat} \Groth(\C)$ be the large set of all discrete Grothendieck fibrations with small fibers. Moreover, recall the notation convention from \cref{subsec:notation} regarding $\Fun(\C^{op},\set)_{/F}$.

\begin{corone}\label{cor:one}
	Let $F: \C^{op} \to \set$ be a functor. Then the isomorphism in \cref{lemma:groth iso} induces a bijection of large sets 
	$$(\Fun(\C^{op},\set)_{/F})^{\sm} \cong \Groth(\C)_{/F}.$$
\end{corone}

\begin{corone} \label{cor:two}
	The isomorphism in \cref{lemma:groth iso} induces a bijection of sets 
	\begin{center}
		\begin{tikzcd}
			 \ds\coprod_{\C \in \cat} \Fun(\C^{op},\set) \arrow[rr, "\int", "\cong"', shift left=1.8] \arrow[dr, "\pi"] & & \Groth \arrow[ll, "\Fib", shift left=1.8] \arrow[dl] \\
			& \cat &
		\end{tikzcd}.
	\end{center} 
\end{corone} 

 We can now move on to construct a functor with value fibrations that avoids the pseudo-functoriality of the pullback. Let $\All: \ssSet^{op} \to \set$ be defined as the composition 
 \begin{equation}\label{eq:all}
 \Fun(\DD^{op}\times\DD^{op},\set)^{op} \xrightarrow{ \ (\int_{\DD\times\DD})^{op} \ } \Groth(\DD\times\DD)^{op} \xrightarrow{ \ \pi^{op} \ } \cat^{op} \xrightarrow{ \ \Fun(-,\set) \ } \set,
 \end{equation} 
 meaning $\All(X) = \Fun(\int_{\DD\times\DD}X,\set)$. Now, we have the following key lemma with regard to $\All$.
 
 \begin{lemone} \label{lemma:bijection}
 	For every bisimplicial set $X$ there is a bijection of sets 
 	$$\Gamma_X: \All(X) \xrightarrow{ \ \cong \ } (\ssSet_{/X})^{\sm}.$$
 \end{lemone}

 \begin{proof}
 	First of all, by restricting the bijection in \cref{cor:two} to the fiber over $\int_{\DD\times\DD}X$, we have the bijection
 	$$\All(X) \cong \Groth(\int_{\DD\times\DD}X).$$
 	Now, there is an evident discrete Grothendieck fibration $\int_{\DD\times\DD}X \to \DD \times\DD$ and so every discrete Grothendieck fibration $\int_{\DD\times\DD}X$ is simply a discrete Grothendieck fibration into $\DD\times\DD$ that factors through $\int_{\DD\times\DD}X$, meaning we have the bijection 
 	$$\Groth(\int_{\DD\times\DD}X) \cong \Groth(\DD\times\DD)_{/\int_{\DD\times\DD}X}.$$
 	Finally, by \cref{cor:one}, we have the bijection 
 	$$\Groth(\DD\times\DD)_{/\int_{\DD\times\DD}X} \cong (\Fun(\DD^{op} \times \DD^{op},\set)_{/X})^{\sm}= (\ssSet_{/X})^{\sm}.$$
 	Combining these three bijections gives us the desired bijection $\All(X) \cong (\ssSet_{/X})^{\sm}$.
 \end{proof}
 
 This bijection is in fact quite well-behaved as we can easily witness by tracing through the definition.
 
 \begin{lemone} \label{lemma:gamma pullback}
 	Let $f: X \to Y$ be a morphism of bisimplicial sets. Then the following square commutes.
 	\begin{center}
 		\begin{tikzcd}
 			\All(Y) \arrow[r, "\Gamma_Y", "\cong"'] \arrow[d, "\All(f)"'] & (\ssSet_{/Y})^{\sm} \arrow[d, "f^*"] \\
 			\All(X) \arrow[r, "\Gamma_X", "\cong"'] & (\ssSet_{/X})^{\sm}
 		\end{tikzcd}.
 	\end{center}
 \end{lemone}

 As a result of this lemma we can think of elements in $\All(X)$ as simplicial spaces over $X$. Now, we want to prove that $\All$ is representable. Define $\ksSet$ as the bisimplicial set obtained by precomposing $\All$ with the Yoneda embedding $\DD^{op} \times\DD^{op} \to \ssSet^{op}$. Concretely, as we can immediately compute $\int_{\DD\times\DD}\Hom_{\DD\times\DD}(-,([n],[l]))= \DD_{/[n]}\times\DD_{/[l]}$, we have $\kssSet(n,l) = \Fun(\DD_{/[n]} \times \DD_{/[l]},\set)$.  We now have the following key result.
 
 \begin{corone} \label{cor:all rep}
 	There is a natural isomorphism $\All \cong \Hom_{\ssSet}(-,\kssSet)$.
 \end{corone}
 
 \begin{proof}
 	The functor $\Hom(-,\kssSet)$ preserves colimits by definition and $\All$ preserves colimits as the three functors in \ref{eq:all} defining $\All$ preserve colimits. Hence, it suffices to observe a natural isomorphism at the level of representables. However, we have 
 	$$\All(\Delta[n,l]) = \Fun(\DD_{/[n]} \times \DD_{/[l]},\set) \cong \Hom_{\ssSet}(\Delta[n,l],\kssSet)$$
 	where the last step follows from the Yoneda lemma, giving us the desired natural isomorphism.
 \end{proof}

The construction of $\All$ and $\kssSet$ is too broad and we often want to restrict it appropriately. We have the following simple observation.
 
 \begin{lemone} \label{lemma:pullback stable}
 	Let $S$ be a (possibly large) set of morphisms in $\ssSet$. The following are equivalent.
 	\begin{enumerate}
 		\item The pullback of a morphism in $S$ (along any morphisms) is in $S$.
 		\item Let $F:\int_{\DD^{op}\times\DD^{op}}X \to \set$ be a functor such that $\Gamma_X(F)$ is in $S$, then for any morphism of bisimplicial sets $f: Y \to X$, $\Gamma_Y(\All(f)(F))$ is in $S$.
 	\end{enumerate}
 \end{lemone}
 
 \begin{proof}
 	Condition $(1)$ corresponds to $f^*:(\ssSet_{/Y})^{\sm} \to (\ssSet_{/X})^{\sm}$ restricting to the full sub set of morphisms in $S$, whereas condition $(2)$ corresponds to $\All(f): \All(Y) \to \All(X)$ restricting similarly. By \cref{lemma:gamma pullback} these two conditions are equivalent as the horizontal morphisms $\Gamma_X,\Gamma_Y$ are bijections. 
 \end{proof}

  We say $S$ is {\it pullback stable} if it satisfies the equivalent conditions in \cref{lemma:pullback stable}. For a given pullback stable set of morphisms $S$, let $\All^S$ be the subfunctor of $\All$ with $F \in \All^S(X)$ if and only if $\Gamma_X(F)$ is in $S$. The functoriality immediately follows from \cref{lemma:pullback stable}. We can similarly define $\kssSet^S$ as the sub-bisimplicial set of $\kssSet$. We now want to deduce a result analogous to \cref{cor:all rep}. For that we need the following additional condition.
  
  \begin{lemone}\label{lemma:local}
  	Let $S$ be a pullback stable set of morphisms in $\ssSet$. The following are equivalent.
  	\begin{enumerate}
  		\item A morphism $Y \to X$ is in $S$ if and only if for all $\Delta[n,l] \to X$, the pullback is in $S$.
  		\item For a functor $F:\int_{\DD^{op}\times\DD^{op}}X \to \set$ we have $\Gamma_X(F)$ in $S$ if and only if every functor $G:\DD_{/[n]} \times \DD_{/[l]} \to \set$ that factors through $F$ we have $\Gamma_{\DD[n,l]}(G)$ is in $S$.
  	\end{enumerate}
  \end{lemone}
  
  \begin{proof}
 	For a given morphism $f: \Delta[n,l] \to X$, \cref{lemma:gamma pullback} gives us the following diagram
 	 	\begin{center}
 		\begin{tikzcd}[column sep=0.5in]
 			\All(\Delta[n,l]) \arrow[r, "\Gamma_{\Delta[n,l]}", "\cong"'] \arrow[d, "\All(f)"'] & (\ssSet_{/\Delta[n,l]})^{\sm} \arrow[d, "f^*"] \\
 			\All(X) \arrow[r, "\Gamma_X", "\cong"'] & (\ssSet_{/X})^{\sm}
 		\end{tikzcd}.
 	\end{center}
 	The assumptions are now direct translations along the bijections $\Gamma_{\Delta[n,l]}$ and $\Gamma_X$, using the fact that $\int_{\DD\times\DD}\Delta[n,l] = \DD_{/[n]} \times \DD_{/[l]}$.
  \end{proof}
  
  A set of morphism that satisfies the equivalent conditions of \cref{lemma:local} is called {\it local}. 
  
  \begin{lemone} \label{lemma:local bijection}
   If $S$ is local, then the bijection given in \cref{cor:all rep} restricts to a bijection $\All^S \cong \Hom_{\ssSet}(-,\kssSet^S)$.	
  \end{lemone}
  
  \begin{proof}
  	As $S$ is local, the sub-functor $\All^S \subseteq \All: \ssSet^{op} \to \set$ is again colimit preserving. The result now follows from the same argument in \cref{cor:all rep}.
  \end{proof}

It useful to have a quick criterion to determine local classes of morphisms with a proof analogous to \cite[Lemma 3.10]{rasekh2017left}.

\begin{corone} \label{cor:lifting local}
	Let $S$ be a set of morphism of bisimplicial sets determined by a right lifting property with respect to a set of morphisms $A \hookrightarrow \Delta[n,l]$. Then $S$ is pullback stable and local. 
\end{corone}

We end this subsection with an elegant example of the previous corollary.

\begin{exone} \label{ex:reedy}
	By \cref{def:reedy}, the large set of Reedy fibrations with small fiber satisfies the condition of \cref{cor:lifting local}. We denote $\All^{\Ree}$ by $\Ree$ and $\kssSet^{\Ree}$ by $\kRee$ and notice that by \cref{lemma:local bijection} we have a natural bijection 
	$$\Ree \cong \Hom_{\ssSet}(-,\kRee).$$
\end{exone} 

\subsection{Minimal Fibrations} \label{subsec:min fib}
In this subsection we introduce minimal Reedy and left fibrations, which play a key role in the construction of strict pullbacks (as discussed in the beginning of \cref{sec:tech}). Recall that a Kan fibration $p: Y \to X$ is a minimal fibration if for any two maps $f,g: \Delta^n \to Y$, such that $f$ is homotopic to $g$ relative to $\partial \Delta^n$ and $pf = pg$, then $f=g$. For more details see \cite[Diagram I.10.1]{goerssjardine2009simplicialhomotopytheory}. We can now generalize this definition directly.

\begin{defone}
	A Reedy fibration of simplicial spaces $Y \to X$ is {\it minimal} if the Kan fibration $Y_n \to X_n \times_{M_nX} M_nY$ is a minimal Kan fibration for all $n\geq0$.
\end{defone}

\begin{remone} \label{rem:min ree local}
	The minimality of a Reedy fibration is determined via equality over a representable $\Delta[n,l]$, meaning a Reedy fibration $p:Y \to X$ is minimal if and only if its pullback along any map $\Delta[n,l] \to X$ is minimal, proving that minimal Reedy fibrations are also pullback stable (\cref{lemma:pullback stable}) and local (\cref{lemma:local}).
\end{remone} 

One key result regarding minimal Kan fibrations is that every Kan fibration $p: Y \to X$ can be restricted to a minimal subfibration $\cMin(p):\cMin(Y) \to X$, such that the inclusion $i: \cMin(Y) \to Y$ has a retract $r:Y \to \cMin(Y)$ over $X$ that is a trivial Kan fibration \cite[Proposition 10.3]{goerssjardine2009simplicialhomotopytheory}. We now have the following generalization.

\begin{propone} \label{prop:factorization}
Every Reedy fibration of bisimplicial sets $p:Y \to X$ can be restricted to a minimal subfibration $\cMin(p):\cMin(Y) \to X$, such that the inclusion $i: \cMin(Y) \to Y$ has a retract $r:Y \to \cMin(Y)$ over $X$ that is a trivial Reedy fibration.
\end{propone}

\begin{proof}
	Using the result for Kan fibrations \cite[Proposition 10.3]{goerssjardine2009simplicialhomotopytheory} inductively, we fix $\cMin(Y)_0 \hookrightarrow Y_0 \twoheadrightarrow \cMin(Y)_0$, using the fact that $Y_0 \to X_0$ is a Kan fibration. Now, assume we have chosen $\cMin(Y)_k \to Y_k \to \cMin(Y)_k$ for $k \leq {n-1}$, which in particular means we can define the matching object $M_n\cMin(Y)$, as it is characterized as a limit of objects $\cMin(Y)_k$ for $k \leq {n-1}$ \cite[Definition 5.2.2]{hovey1999modelcategories}. We now have the following diagram (here arrows labeled with $min$ are minimal fibrations)
	\begin{center}
		\begin{tikzcd}[row sep=0.5in, column sep=0.5in]
			\cMin(Y)_n \arrow[r, hookrightarrow] \arrow[rrr, bend left=17, equal]& Y_n \arrow[rr, bend left=15, twoheadrightarrow, "\simeq" description] \arrow[r, twoheadrightarrow, "p_n" description] \arrow[dr, twoheadrightarrow] & M_nY \times_{M_nX} \cMin(Y)_n \arrow[r, twoheadrightarrow, "\simeq" description] \arrow[d, twoheadrightarrow, "min" description] \arrow[dr, phantom, "\ulcorner", very near start] & \cMin(Y)_n \arrow[d, twoheadrightarrow, "min" description] & \\
			& & M_nY \times_{M_nX} X_n \arrow[r, twoheadrightarrow, "\simeq" description] \arrow[d] \arrow[dr, phantom, "\ulcorner", very near start] & M_n\cMin(Y) \times_{M_nX}  X_n \arrow[r, twoheadrightarrow] \arrow[d] \arrow[dr, phantom, "\ulcorner", very near start] & X_n \arrow[d]  \\
			& & M_nY \arrow[r, twoheadrightarrow, "\simeq" description] & M_n\cMin(Y) \arrow[r, twoheadrightarrow] & M_nX
		\end{tikzcd}.
	\end{center}
	Here $\cMin(Y)_n$ is chosen as a factorization of the fibration $Y_n \twoheadrightarrow M_n\cMin(Y) \times_{M_nX}X_n$ (using  \cite[Proposition 10.3]{goerssjardine2009simplicialhomotopytheory}), meaning $\cMin(Y)_n \to  M_n\cMin(Y) \times_{M_nX}X_n$ is a minimal fibration, and notice the map $p_n$ is in fact a Kan fibration as it obtained as the limit of the diagram
	\begin{center}
		\begin{tikzcd}
			Y_n \arrow[d, twoheadrightarrow] \arrow[r] & Y_n \arrow[d, twoheadrightarrow] & Y_n \arrow[l] \arrow[d, twoheadrightarrow] \\
			M_nY \times_{M_nX} X_n \arrow[r] & M_n\cMin(Y) \times_{M_nX} X_n & M_n \arrow[l]
		\end{tikzcd}
	\end{center}
	and is in fact a trivial fibration by $2$-out-of-$3$. Hence, we have obtained $\cMin(Y)_n \hookrightarrow Y_n$ with all desired conditions, which completes our induction step.
\end{proof}

Our construction works well for an individual Reedy fibration, however, we would like to have a construction of the minimal Reedy fibration that is consistent with pullback. This requires us to make a globally consistent choice of minimal fibrations for all Reedy fibrations at once, which we achieve in the following way.

Let $\Reemin \subseteq \Ree$ be the subset of minimal Reedy fibration with small fiber and denote the corresponding subobject of the bisimplicial set of $\kRee$ by $\kReemin$. By \cref{rem:min ree local} minimal Reedy fibrations are local and so we have a bijection
\begin{equation} \label{eq:min Ree}
 \Reemin \cong \Hom_{\ssSet}(-,\kReemin)
\end{equation}
and in particular we have a bijection of sets $\Reemin(\kRee) \cong \Hom_{\ssSet}(\kRee,\kReemin)$. Now, applying \cref{prop:factorization} to the Reedy fibration over $\kRee$ that corresponds to the identity map in \cref{ex:reedy} we obtain a minimal Reedy fibration over $\kRee$ that by \ref{eq:min Ree} corresponds to a map of bisimplicial sets 
\begin{equation} \label{eq:min ree}
 \kMin:\kRee \to \kReemin.
\end{equation}
 Let $\cMin: \Ree \to \Reemin$ be the map represented by $\kMin$. This map sends every Reedy fibration to its corresponding minimal Reedy fibration constructed in \cref{prop:factorization}. Indeed, for every Reedy fibration $Y \to X$, we have the following diagram, where $\kRee_* \to \kRee$ is the Reedy fibration corresponding to the identity morphism in \cref{ex:reedy}
\begin{center}
	\begin{tikzcd}
		Y \arrow[r] \arrow[d, twoheadrightarrow, "\simeq"] \arrow[dr, phantom, "\ulcorner", very near start] & \kRee_* \arrow[d, twoheadrightarrow, "\simeq"] \\
		\cMin(Y) \arrow[r]  \arrow[d, twoheadrightarrow, "min"] \arrow[dr, phantom, "\ulcorner", very near start] & \cMin(\kRee_*) \arrow[d, twoheadrightarrow, "min"]  \\
		X \arrow[r] & \kRee 
	\end{tikzcd}, 
\end{center} 
where we are using the fact that both minimal Reedy fibrations and trivial fibrations are stable under pullback.

How does the minimality construction interact with left fibrations?

\begin{lemone} \label{lemma:min left}
	Let $L \to X$ be a left fibration, then $\cMin(L) \to X$ is also a left fibration.
\end{lemone} 

\begin{proof}
 We already know that $\cMin(L) \to X$ is a Reedy fibration, hence, by \cref{def:left fib} it suffices to observe that $\cMin(L)_n \to \cMin(L)_0 \times_{X_0} X_n$ is an equivalence, which follows directly from the fact that $L$ is a left fibration and $L_n \simeq \cMin(L)_n$ for all $n \geq 0$ (\cref{prop:factorization}).
\end{proof}

Finally, the key concept that makes minimal Kan fibrations so useful is that every equivalence between two minimal Kan fibrations is in fact an isomorphism \cite[Lemma 10.4]{goerssjardine2009simplicialhomotopytheory} and we have the following analogous result.

\begin{lemone} \label{lemma:reedy equiv strict}
	Let $p:Y \to X, q:Z \to X$ be two minimal Reedy fibrations and $f: Y \to Z$ a map over $X$. Then $f$ is a Reedy equivalence if and only if it is an isomorphism. 
\end{lemone}

\begin{proof}
	For all $n\geq0$ we have the following commutative diagram 
	\begin{equation} \label{eq:square min fib}
		\begin{tikzcd}
			Y_n \arrow[r, "f_n"] \arrow[d, twoheadrightarrow] & Z_n \arrow[d, twoheadrightarrow] \\
			M_nY \times_{M_nX} X_n \arrow[r] & M_nZ \times_{M_nX} X_n 
		\end{tikzcd},
	\end{equation}
	where the two vertical morphisms are minimal Kan fibrations. We now prove the statement by induction. $f_0$ is an isomorphism by \cite[Lemma 10.4]{goerssjardine2009simplicialhomotopytheory}. Now assuming it holds for all $k < n$ it follows that the bottom map in \ref{eq:square min fib} is an isomorphism and hence \cite[Lemma 10.4]{goerssjardine2009simplicialhomotopytheory} implies again that $f_n$ is an isomorphism.
\end{proof}

\section{The Complete Segal Space of Spaces} \label{sec:spaces}
In this section we finally make the intuition outlined in \cref{subsec:higher cat fib} precise and construct the desired complete Segal space of spaces using left fibration. Let $\LFib$ be the large set of left fibrations with small fibers. By \cref{def:left fib} and \cref{cor:lifting local}, $\LFib$ is local (\cref{lemma:local}) and so we can take the sub-functor of $\All:\ssSet^{op} \to \set$ with value left fibrations that we denote by $\LFib(-):\ssSet^{op} \to \set$, which, by \cref{lemma:local bijection}, is represented by a bisimplicial set that we denote by $\kS$, meaning we have a natural bijection 
\begin{equation} \label{eq:lfib}
 \LFib(-) \cong \Hom_{\ssSet}(-,\kS).
\end{equation}
We now want to prove that $\kS$ is a complete Segal space of spaces. Before we can get to the main result we need the appropriate  lemma that helps us understand extension properties of trivial fibrations. For that we can directly generalize \cite[Lemma 5.1.20]{cisinski2019highercategories} to bisimplicial sets.

\begin{lemone} \label{lemma:triv fib}
	Let $p: Y \to A$ be a trivial Reedy fibration of bisimplicial sets and $i:A \to B$ an inclusion of bisimplicial sets. Then $i_*p$ is a trivial Reedy fibration and $p \xrightarrow{ \ \cong \ } i^*i_*p$.
\end{lemone}

\begin{propone} \label{prop:spaces}
	The bisimplicial set $\kS$ is a complete Segal space.
\end{propone}

\begin{proof}
	The argument is analogous to \cite[Theorem 2.2.1]{kapulkinlumsdaine2021kanunivalent} and \cite[Theorem 5.2.10]{cisinski2019highercategories}. By \cref{rem:css lift}, we need to prove that for every trivial cofibration $i:A \to B$ in the complete Segal space model structure, the induced map  $\Hom_{\ssSet}(B,\kS) \to \Hom_{\ssSet}(A,\kS)$ is surjective. By \ref{eq:lfib}, this is equivalent to $\LFib(B) \to \LFib(A)$ being surjective, which concretely means proving that every left fibration over $A$ is the pullback of a left fibration over $B$ via $i$.
	
	Fix a left fibration $p: L \to A$. We now have the following diagram 
	\begin{center}
		\begin{tikzcd}[row sep=0.5in, column sep=0.7in]
			L \arrow[dd, "p"', bend right=40, twoheadrightarrow] \arrow[d, twoheadrightarrow, "\simeq", "r"'] \arrow[rr, hookrightarrow] \arrow[dr, phantom, "\ulcorner", very near start] & & (\hat{r}j)_*L \arrow[d, twoheadrightarrow, "\simeq"', "(\hat{r}j)_*r"] \\
			\cMin(L) \arrow[d, twoheadrightarrow, "\cMin(p)"] \arrow[r, "j", "\simeq"', hookrightarrow] \arrow[dr, phantom, "\ulcorner", very near start] & \hat{L} \arrow[d, twoheadrightarrow, "\hat{p}"] \arrow[r, twoheadrightarrow, "\simeq"', "\hat{r}"] & \cMin(\hat{L}) \arrow[dl, twoheadrightarrow, "\cMin(\hat{p})"] \\  
			A \arrow[r, "i", hookrightarrow] & B
		\end{tikzcd},
	\end{center}
	where $\cMin(p)$ is a minimal left fibration (\cref{lemma:min left}) and $r$ a trivial fibration (\cref{prop:factorization}), $\hat{p}$ is a left fibration and $j$ a trivial complete Segal space equivalence (\ref{eq:derived unit}), $\hat{r}$ a trivial fibration, $\cMin\hat{p}$ a minimal left fibration (again by \cref{prop:factorization}) and $(\hat{r}j)_*r$ is a trivial fibration (\cref{lemma:triv fib}). Now, by the properties of the homotopy pullback square \ref{eq:derived unit}, the map $\cMin(L) \to i^*\cMin(\hat{L})$ induced by the pullback is a Reedy equivalence and hence, by \cref{lemma:reedy equiv strict}, a bijection and, by \cref{lemma:triv fib}, the top rectangle is a pullback. Hence $p$ is the pullback of the left fibration $\cMin(\hat{p}) \circ (\hat{r}j)_*r:(\hat{r}j)_*L \to B$ and we are done.
\end{proof} 

\begin{remone} \label{rem:kazhdan}
	There is a similar result in \cite[Theorem 2.2.11]{kazhdanvarshvsky2014yoneda} without addressing the functoriality of the construction given that their definition of the simplicial space uses pullbacks \cite[Main construction 2.2.3]{kazhdanvarshvsky2014yoneda}.
\end{remone}

We now want to understand the mapping spaces of $\kS$. For that we need the following strictification. Combining the locality of left fibrations (\ref{eq:lfib}) and minimal Reedy fibrations (\ref{eq:min Ree}) it follows, by \cref{lemma:local bijection}, that the set of minimal left fibrations is local and so we have a sub bisimplicial set of $\kS$, that we denote by $\kS^{min} \hookrightarrow \kS$, and natural bijection 
\begin{equation} \label{eq:bij min left}
	\LFibmin(-) \cong \Hom_{\ssSet}(-,\kS^{min}).
\end{equation}
On the other hand the map $\kMin:\kRee \to \kRee^{min}$ defined in \ref{eq:min ree}, by \cref{lemma:min left}, restricts to a map $\kMin:\kS \to \kS^{min}$. We now have the following result with regard to these two maps.

\begin{lemone} \label{lemma:s vs s min}
	The maps $\kS^{min} \to \kS \to \kS^{min}$ are equivalences of complete Segal spaces. 
\end{lemone} 

\begin{proof}
	First of all  the composition $\kS^{min} \to \kS \to \kS^{min}$ is the identity as it takes every minimal left fibration to itself. Hence, $\kS^{min}$ is a retract of $\kS$ and so a complete Segal space, Next, we prove that $\kS \to \kS^{min}$ is a trivial fibration and this implies that the inclusion is an equivalences as well. 
	Following \cref{def:reedy}, we need to prove for every inclusion of simplicial spaces $i:A \to B$ the following diagram has a lift
	\begin{center}
		\begin{tikzcd}
			A \arrow[r] \arrow[d, "i"] \arrow[dr, phantom, "\ulcorner", very near start] & \kS \arrow[d, "\kMin"] \\
			B \arrow[r] & \kS^{min}
		\end{tikzcd},
	\end{center}
	which, by \ref{eq:lfib} and \ref{eq:bij min left}, is equivalent to the map
	$$\LFib(B) \to \LFib(A) \times_{\LFibmin(A)} \LFibmin(B)$$
    being surjective. Unwinding the definitions this means we have the data of the following diagram 
    \begin{center}
    	\begin{tikzcd}
    		L \arrow[d, twoheadrightarrow, "r", "\simeq"'] & \\
    		\cMin(L) \arrow[r, "\hat{p}^*i"] \arrow[d, twoheadrightarrow, "p"] \arrow[dr, phantom, "\ulcorner", very near start] & \hat{L} \arrow[d, twoheadrightarrow, "\hat{p}"] \\
    		A \arrow[r, "i"] & B
    	\end{tikzcd},
    \end{center}
	where $r$ is a trivial fibration and $p,\hat{p}$ are minimal fibrations and we need to find a left fibration $\tilde{p}:L \to B$, such that $i^*\tilde{p} = pr$ and $\cMin(\tilde{p}) = \hat{p}$. However, by \cref{lemma:triv fib}, this is given by $(\hat{p}^*i)_*r:(\hat{p}^*i)_*L \to \hat{L}$.
\end{proof}

We now want to use \cref{lemma:s vs s min} to better understand the mapping spaces (\ref{eq:map}) of $\kS$. This requires us to better understand minimal left fibrations over $\Delta[1,0]$. First we introduce a notation that will be useful in the next proofs.

\begin{notone} \label{not:arrow}
	Let $X \to \Delta[1,l]$ be a map of bisimplicial sets. We use the following three notational conventions
	\begin{itemize}
		\item $X_{/0}=\Map_{/\Delta[1,l]}(d^1:\Delta[0,l]\to \Delta[1,l],X)$
		\item $X_{/1}=\Map_{/\Delta[1,l]}(d^0:\Delta[0,l]\to \Delta[1,l],X)$
		\item $X_{/01}=\Map_{/\Delta[1,l]}(\id,X)$
	\end{itemize}
	Notice, the two maps $d^0,d^1:\Delta[0,l]\to \Delta[1,l]$ induce maps of simplicial sets $s:X_{/01} \to X_{/0}$ and $t: X_{/01} \to  _{/1}$.
\end{notone}

\begin{lemone} \label{lemma:lfib min arrow}
	Let $L$ and $L' $ be two left fibrations over $\Delta[1,l]$ such that the two morphisms $t_L:L_{/01} \to L_{/1}$, $t_{L'}:L_{/01}' \to L_{/1}'$ are Kan equivalent morphisms. Then $\cMin(L) = \cMin(L')$.
\end{lemone}

\begin{proof}
 By \cref{lemma:reedy equiv strict} the result will follow if we can prove that $L$ and $L'$ are equivalent over $\Delta[1,l]$. Now, the projection map $\Delta[1,l] \to \Delta[1,0]$ is a level-wise equivalence and so we, can without loss of generality assume that $l=0$. Denote the equivalence $t_L \to t_{L'}$ by $\alpha$.
 
 We now use the adjunction $(\sint_{[1]},\sH_{[1]})$ as defined in \cite[Lemma 4.9]{rasekh2017left}. Concretely, by definition of $\sH_{[1]}$ we have
 $\sH_{[1]}(L) = t_L$, $\sH_{[1]}(L') = t_{L'}$ and by \cite[Theorem 4.18]{rasekh2017left} the counit of the adjunction gives us a level-wise equivalence and so we have
 $$L \simeq \sint_{[1]}\sH_{[1]} L = \sint_{[1]} t_L\xrightarrow[\simeq]{\sint_{[1]}\alpha} \sint_{[1]} t_{L'} = \sint_{[1]}\sH_{[1]} L' \simeq L' $$
 and hence we are done.
\end{proof}

 As explained in \cref{subsec:nerves} one easy way to construct an $(\infty,1)$-category of spaces is via nerves of Kan enriched categories. We hence want to compare $\kS$ to the nerve of the Kan enriched category of Kan complexes. Following the notational convention of \cref{subsec:notation}, let  $\uKan$ be the Kan enriched category of Kan complexes, which has objects Kan complexes and for two objects $K,L$ we have 
\begin{equation}\label{eq:map kan}
	\Map_{\uKan}(K,L)_n = \Hom_{\uKan}(K \times \Delta[n],L) \cong \Hom_{/\Delta[l]}(K \times \Delta[n],L \times \Delta[n])
\end{equation}
For further details regarding this Kan enriched category see \cite[Subsection I.5]{goerssjardine2009simplicialhomotopytheory}.

\begin{remone} \label{rem:kanmin}
  Analogous to \cref{lemma:s vs s min} we say a Kan complex $K$ is {\it minimal} if $K =\cMin(Y)_0$ (\ref{eq:min ree}) for some Reedy fibrant simplicial space $Y$ and define $\uKanmin$ as the full Kan enriched subcategory of $\uKan$ with objects the minimal Kan complexes and notice the inclusion $\uKanmin \to \uKan$ is an equivalence of Kan enriched categories, in the sense of \cite{bergner2007bergnermodelcat}, as every Kan complex is equivalent to a minimal Kan complex (\cref{prop:factorization}).
\end{remone}

We now want to apply the nerve to these Kan enriched categories. For a given simplicially enriched category $\C$, let $N_\Delta\C$ be the bisimplicial set with $N_\Delta\C_0 = \Obj_\C$ and 
\begin{equation} \label{eq:nerve}
 N_\Delta\C_n = \coprod_{X_0,...,X_n} \Map_\C(X_0,X_1) \times ... \times \Map_\C(X_{n-1},X_n).
\end{equation} 

While $N_\Delta\C$ is not a complete Segal space it is in fact a Segal category \cite[Proposition 8.3]{bergner2007threemodels} (where the nerve $N_\Delta$ is denoted $R$ instead) and we can characterize their equivalences via {\it Dwyer-Kan} equivalences meaning we have the following special instance of \cite[Theorem 7.1]{bergner2007threemodels}.

\begin{propone} \label{prop:segal cat}
	Let $\C$ be a Kan enriched category and $W$ a complete Segal space. A map $F:N_\Delta\C \to W$ is an equivalence in the complete Segal space model structure if $\Obj_\C \to W_{00}$ is surjective and for objects $x,y$ in $\C$, the induced map $\Map_\C(x,y) \to \Map_W(Fx,Fy)$ is a Kan equivalence.
\end{propone} 

We now have the necessary pieces to prove the main result.

\begin{theone} \label{the:main spaces}
	There is an equivalence between the complete Segal space $\kS$ and the strict Segal category $N_\Delta\uKan$. Moreover, we have a natural bijection $\LFib \cong \Hom(-,\kS)$.
\end{theone}
	
\begin{proof}
	 We have established that $\kS$ is a complete Segal space in \cref{prop:spaces} and the bijection in \ref{eq:lfib} and so we only need to prove the equivalence with $N_\Delta\uKan$. Let $\sL \to N_\Delta\uKan$ be the map of bisimplicial sets with 
	 \begin{equation} \label{eq:special lfib}
      	\sL_n = \coprod_{X_0,...,X_n} X_0 \times \Map(X_0,X_1) \times \cdots \times \Map(X_{n-1},X_n)	 
     \end{equation}  
	 with boundary maps given by composition and projection. The evident projection map $\pi_1:\sL \to N_\Delta\uKan$ gives us the strict pullback of simplicial set $\sL_n \cong \sL_0 \times_{\Obj_{\uKan}} N_\Delta\uKan_n$ and so, by \cref{rem:strict left}, the Reedy fibrant replacement $\pi_1:\hat{\sL} \to N_\Delta\uKan$ is a left fibration. By \ref{eq:lfib}, this induces a functor $N_\Delta\uKan \to \kS$. We want to prove it is a complete Segal space equivalence. 
	 
	 By \cref{lemma:s vs s min}, it suffices to prove that the composition $N_\Delta\uKan \to \kS^{min}$ is an equivalence. Moreover, by \cref{rem:kanmin}, we can further reduce it to showing that $\Lift:N_\Delta\uKanmin \to \kS^{min}$ is a Dwyer-Kan equivalence. By \cref{rem:lfib constant} and \cref{lemma:reedy equiv strict} every left fibration over $L \to \Delta[0,0]$ is uniquely determined by the minimal Kan complex $L_0$ and so $\Lift$ is surjective on object and so, by \cref{prop:segal cat}, it suffices to prove $\Lift$ induces an equivalence of mapping spaces. We will in fact show $\Lift_1$ induces a bijection on mapping spaces.
	 
	 Before we proceed, we will understand this map better. By construction, $\Lift$ corresponds under the bijection \ref{eq:bij min left} to the minimal left fibration $\hat{\sL}^{min}=\cMin(\hat{\sL}) \times_{\uKan} \uKanmin$ over $\uKanmin$, meaning $\Lift(\sigma)= \sigma^*\hat{\sL}^{min} \to \Delta[n,l]$, for every $\sigma:\Delta[n,l] \to N_\Delta\uKanmin$. By definition of the simplicial nerve (\ref{eq:nerve}) $(N_\Delta\uKanmin)_{1l}$ is the set of $l$-morphisms in the simplicially enriched category $\uKanmin$ and, by the explanation above, $\Lift$ will take $f: X \times \Delta[l] \to Y \times \Delta[l]$ over $\Delta[l]$ to the minimal left fibration $\Lift(f) \to \Delta[1,l]$ with the following characteristics (using \cref{not:arrow}):
	 \begin{itemize}
	 	\item $\Lift(f)_{/0} = X, \Lift(f)_{/1}= Y$,
	 	\item $\Lift(f)_{/01}$ is the minimal Reedy fibration over $X \times Y$ equivalent to the map $(\id,f):X \to X\times Y$, which determines $\Lift(f)_{/01}$ uniquely.
	 \end{itemize}
	  The equivalence $X \simeq \Lift(f)_{/01}$ in particular implies that 
	  \begin{equation}\label{eq:min}
	  	\cMin(\Lift(f)_{/01}) = X
	  \end{equation}
  	  and $X =\cMin(\Lift(f)_{/01}) \to X \times Y \xrightarrow{ \ \pi_2 \ } Y$ is given by $f$.
	 
	 We will now prove $\Lift_1$ is a bijection by constructing an inverse $\Sec: (\kS^{min})_1 \to N_\Delta\uKanmin_1$. For a given minimal left fibration $L \to \Delta[1,l]$, let $\Sec(L)= \cMin(L_{/01}) \to L_{/1}$ over $\Delta[l]$. We now prove that $\Sec$ is injective and that $\Sec$ is a left inverse of $\Lift_1$ which will prove they are inverses.
	 
	 The statement of \cref{lemma:lfib min arrow} and the definition of $\Sec(L)$ as $L_{/01} \to L_{/1}$ immediately implies that $\Sec$ is injective on minimal left fibrations. Moreover, we observed in \ref{eq:min} that for a given morphism of minimal left fibrations $f: X \times\Delta[l] \to Y \times\Delta[l]$ over $\Delta[l]$, $X= \cMin(\Lift(f)_{/01})$ and the composition 
	 $$X=\cMin(\Lift(f)_{/01}) \hookrightarrow \Lift(f)_{/01} \to \Lift(f)_{/1}= Y$$ 
	 is given by $f$. This proves that $\Sec\circ \Lift_1$ is the identity and hence we are done.
\end{proof} 

\begin{remone} \label{rem:rfib fib}
 We constructed a complete Segal space of spaces using left fibrations. However, based on \cref{rem:rfib}, we could have also used right fibrations to construct a bisimplicial set $\kS^{R}$ with $\kS^{R}_{nl}$ given by right fibrations over $\Delta[n,l]$. Now, \cref{rem:rfib} implies that $(-)^{op}$ induces a bijection between $\kS$ and $\kS^R$ that flips the directionality of the morphisms, immediately implying that $\kS^{R}$ is just the opposite complete Segal space of spaces, $\kS^{op}$. 
\end{remone}

\begin{remone} \label{rem:kapulkin}
	One implication of \cref{the:main spaces} is that $\kS$ is complete, which in particular implies the space of equivalences from $X$ to $Y$ is equivalent to the space of left fibrations over $\Delta[0,\bullet]$ from $X$ to $Y$. As left fibrations $L \to \Delta[0,n]$ are uniquely determined by Kan fibrations over $\Delta[n]$ \cite[Theorem 3.17]{rasekh2017left}, this means we are getting an equivalence between the space of equivalences and the space of Kan fibrations with specific fibers, meaning we get an alternative proof to the simplicial univalence of the universe of Kan fibrations \cite[Theorem 3.4.1]{kapulkinlumsdaine2021kanunivalent}.
\end{remone}

\section{The Complete Segal Space of (Complete) Segal Spaces}
We now want to generalize the results from the previous section from a complete Segal space of spaces to a complete Segal space of complete Segal spaces. Here we rely on the complete Segal object approach to Cartesian fibrations as outlined in the beginning of \cite{rasekh2021cartfibmarkedvscso}. The benefit of this approach is that we can use the results of the previous sections level-wise to immediately deduce the desired results. 

Let us start with the appropriate generalizations of \cref{subsec:simplicial}. The category of {\it trisimplicial sets} is denoted by $\sssSet$ and the generators are denoted by $\Delta[-,-,-]$. We can, analogous to \cref{def:reedy}, characterize Reedy fibrations via right lifting properties against certain morphisms with codomain $\Delta[k,n,l]$. Again, Reedy fibrations are part of a model structure with cofibrations given by inclusions of trisimplicial sets and equivalences given by level-wise Kan equivalences. In particular all trivial Reedy fibrations are Reedy weak equivalences. See \cite[Subsection 1.8]{rasekh2021cartfibcss} for more details. 

\begin{notone} \label{not:trisimplicial}
 For a trisimplicial set $X$, we use $X_k$ to denote the bisimplicial set with $(X_k)_{nl} = X_{knl}$ and use $X_{kn}$ to denote the simplicial set with $(X_{kn})_l = X_{knl}$.	 
\end{notone} 

We have analogous generalizations of left fibrations (\cref{def:left fib}).
\begin{defone} \label{def:standard embedding}
	Let $\iota_{st}: \ssSet \to \sssSet$, the {\it standard embedding}, take a bisimplicial set $X$ to the trisimplicial set $\iota_{st}(X)$ characterized as $\iota_{st}(X)_{knl} = X_{nl}$. Moreover, for a map of trisimplicial sets $p:Y \to \iota_{st}X$, we denote by $p_k:Y_k \to X$ the map of bisimplicial sets defined as $(p_k)_{nl} = p_{knl}$.
\end{defone}

\begin{defone} \label{def:reelfib}
	Let $X$ be a bisimplicial set and $p:L \to X$ be a Reedy fibration of trisimplicial sets. 
	\begin{itemize}
		\item $p$ is a {\it Reedy left fibration} if $p_k:L_k \to X$ (\cref{not:trisimplicial}) is a left fibration for all $k \geq 0$.
		\item $p$ is a {\it Segal coCartesian fibration} if it is a Reedy left fibration and it satisfies the {\it Segal condition}, meaning the map 
		$$L_{kn} \to L_{1n} \times_{L_{0n}} ... \times_{L_{0n}} L_{1n}$$
		is a Kan equivalence for all $k\geq 0$ and  $n \geq 2$.
		\item $p$ is a {\it coCartesian fibration} if it is a Segal coCartesian fibration and it satisfies the {\it completeness condition}, meaning the map
		$$L_{k0} \to L_{k3} \times_{L_{k1} \times L_{k1}} L_{k0} \times L_{k0}$$
		is a Kan equivalence for all $k\geq 0$.
	\end{itemize}
\end{defone}

The Segal condition and completeness condition are evident analogues to the ones used to define complete Segal spaces (\cref{def:css}) and so we would expect a close connection. Indeed we have the following remark.

\begin{remone} \label{rem:fiberwise cart fib}
 By \cite[Theorem 4.7]{rasekh2021cartfibcss}, if $p: L \to X$ is a Reedy left fibration, then $p$ is Segal coCartesian fibration if and only if for every $x:\Delta[0,0] \to X$, the fiber $x^*L_0$ is a Segal space and similarly between coCartesian fibrations and complete Segal spaces. 	
\end{remone}

\begin{remone} \label{rem:lfib constant s}
	Following \cref{rem:lfib constant} a Reedy left fibration over $\Delta[0,0]$ is a homotopically constant trisimplicial set. Similarly, a (Segal) coCartesian fibration over $\Delta[0,0]$ is a homotopically constant trisimplicial set $L$ such that $L_0$ is a (complete) Segal space.
\end{remone}
 
We have a similar interaction between Reedy left fibrations and complete Segal equivalences as in \ref{eq:derived unit}, meaning every Reedy left fibration $L \to A$ can be obtained as a homotopy pullback square of a Reedy left fibration $\hat{L} \to B$
\begin{equation} \label{eq:derived unit s}
	\begin{tikzcd}
		L \arrow[r] \arrow[d, twoheadrightarrow] \arrow[dr, phantom, "\ulcorner", very near start] & \hat{L} \arrow[d, twoheadrightarrow] \\
		A \arrow[r, "i"] & B
	\end{tikzcd}
\end{equation}
where $A \to B$ is a trivial complete Segal space cofibration \cite[Theorem 4.12]{rasekh2021cartfibcss}

\begin{remone}\label{rem:rfib s}
  Generalizing \cref{rem:rfib}, we can analogous to \cref{def:reelfib} define {\it Reedy right fibrations} as Reedy fibrations $p:R \to X$, such that for all $k$, $p_k:R_k \to X$ is a right fibration. Moreover, we then define {\it (Segal) Cartesian fibrations} as Reedy right fibrations that satisfy the Segal and completeness condition, as described in \cref{def:reelfib}.
  Moreover, by \cref{rem:rfib}, a map $R \to X$ is a Reedy right, Segal Cartesian or Cartesian fibration if and only if $R^{op} \to X^{op}$ is a Reedy left, Segal coCartesian or coCartesian fibration, respectively. Here $(-)^{op}:\sssSet \to \sssSet$ takes $(-)^{op}$ defined in \cref{rem:rfib} level-wise, meaning it is defined as $\Fun(\DD^{op},(-)^{op})$. 
\end{remone}
We move on to appropriate generalizations of \cref{subsec:fib functor}. Let $\All:\sssSet^{op} \to \set$ be the functor $\All(X) = \Fun(\int_{\DD\times\DD\times\DD}X,\set)$. We can restrict this functor by precomposing with the standard embedding $(\iota_{st})^{op}:\ssSet^{op} \to \sssSet^{op}$ to define $\sAll:\ssSet^{op} \to \set$. Similar to \cref{lemma:bijection} for every bisimplicial set $X$ we have a bijection 
$$\sAll(X) \cong (\sssSet_{/\iota_{st}X})^{\sm}. $$
Now, define $\ksssSet: \DD^{op} \times \DD^{op} \to \set$ as $\sAll$ precomposed with the Yoneda embedding and, analogous to \cref{cor:all rep} we have a natural isomorphism 
\begin{equation} \label{eq:equiv s}
\sAll(-) \cong \Hom_{\ssSet}(-,\ksssSet).
\end{equation} 
Finally, if $S$ is a local set of morphisms, then, similar to \cref{lemma:local bijection}, this bijection restricts to a bijection 
\begin{equation} \label{eq:local equiv}
 \sAll^S(-) \cong \Hom_{\ssSet}(-,\ksssSet),
\end{equation} 
where $\sAll^S(X) \subseteq \sAll(X)$ is the sub-functor of objects over $\iota_{st}X$ that are in $S$ and $\ksssSet^S$ is again the restriction of $\sAll^S$ along the Yoneda embedding. Finally, we have the analogue of \cref{cor:lifting local} for trisimplicial sets.

\begin{corone} \label{cor:lifting local s}
	Let $S$ be a set of morphism of trisimplicial sets determined by a right lifting property with respect to a set of morphisms $A \hookrightarrow \Delta[k,n,l]$. Then $S$ is pullback stable and local. 
\end{corone}

We can use the corollary to generalize \ref{eq:lfib}. By \cref{def:reelfib}, Reedy left fibrations are local and so we get a bisimplicial set $\ksssSet^{\ReeLFib}$ that we denote by $\ksS$ and a natural bijection 
\begin{equation} \label{eq:reelfib}
\ReeLFib(-) \cong \Hom_{\ssSet}(-,\ksS).
\end{equation}
We want to prove that $\ksS$ is a complete Segal space of Reedy fibrant simplicial spaces. This requires us to understand minimal Reedy fibration and minimal Reedy left fibrations of trisimplicial sets, similar to \cref{subsec:min fib}. A Reedy fibration of trisimplicial sets $Y \to X$ is {\it minimal} if for all $k,n \geq 0$ the induced map of simplicial sets (using \cref{not:trisimplicial}) 
$$Y_{kn} \to M_{kn}Y \times_{M_{kn}X} X_{kn}$$
is a minimal Kan fibration. Now using the same inductive argument used in \cref{prop:factorization} (now on two variables $k,n \geq 0$), we have the following result.

\begin{propone} \label{prop:factorization s}
	Every Reedy fibration of trisimplicial sets $p:Y \to X$ can be restricted to a minimal subfibration $\cMin(p):\cMin(Y) \to X$, such that the inclusion $i: \cMin(Y) \to Y$ has a retract $r:Y \to \cMin(Y)$ over $X$ that is a trivial Reedy fibration.
\end{propone}

Applying \cref{prop:factorization s} to  \ref{eq:min ree} we obtain a natural morphism 
\begin{equation} \label{eq:min ree s}
	\scMin: \Ree \to \Reemin
\end{equation}
that naturally takes each Reedy fibration of trisimplicial sets to a minimal Reedy fibration.

It follows from \cref{lemma:min left} that if $L \to X$ is a Reedy left fibration, then $\cMin(L) \to X$ is a Reedy left fibration. We can now use this result to study the bisimplicial set $\ksS$. Before that we need the following last lemma, which follows analogous to \cref{lemma:triv fib}.

\begin{lemone} \label{lemma:triv fib s}
	Let $p: Y \to A$ be a trivial Reedy fibration of trisimplicial sets and $i:A \to B$ an inclusion of trisimplicial sets. Then $i_*p$ is a trivial Reedy fibration and $p \xrightarrow{ \ \cong \ } i^*i_*p$.
\end{lemone}

\begin{propone} \label{prop:spaces s}
	The bisimplicial set $\ksS$ is a complete Segal space.
\end{propone}

\begin{proof}
	We will follow the steps given in \cref{prop:spaces}. The bijection \ref{eq:reelfib} allows us to again reduce the proof to showing that every Reedy left fibration $L \to A$ can be obtained as the pullback of a Reedy left fibration $R \to B$ along a trivial complete Segal space cofibration $i:A \to B$. We can obtain this lift using the same diagram as in the proof of \cref{prop:spaces} this time relying on \cref{prop:factorization s} whenever we need a minimal Reedy left fibration, \ref{eq:derived unit s} when we need to extend $L$ along $i:A \to B$ and \cref{lemma:triv fib s} when we need a trivial Reedy fibration. 
\end{proof}

 Now analogous to \ref{eq:bij min left}, we get a complete Segal space $\ksS^{min} \to \ksS$ consisting of minimal Reedy left fibrations over $\Delta[n,l]$ with a retract $\ksS \to \ksS^{min}$ along with a bijection 
\begin{equation} \label{eq:lfib min}
	\ReeLFib^{min}(-) \cong \Hom_{\ssSet}(-,\ksS^{min}).
\end{equation}
Using this bijection along with \cref{lemma:triv fib s} in the proof of \cref{lemma:s vs s min} we obtain the following result.

\begin{lemone} \label{lemma:s vs s min s}
	The maps $\ksS^{min} \to \ksS \to \ksS^{min}$ are equivalences of complete Segal spaces. 
\end{lemone} 

We now want to use these result to finally prove that $\ksS$ is in fact the complete Segal space of Reedy fibrant bisimplicial sets. Let $\uRee$ denote the Kan enriched category of Reedy fibrant simplicial spaces, where, by analogy with \ref{eq:map kan}, the mapping spaces for two Reedy fibrant simplicial spaces are defined as follows:
\begin{equation}\label{eq:map reedy}
	\Map_{\uRee}(K,L)_n = \Hom_{\uRee}(K \times \Delta[0,n],L) \cong \Hom_{/\Delta[0,n]}(K \times \Delta[0,n],L \times \Delta[0,n]).
\end{equation}

\begin{remone} \label{rem:reemin}
	Similar to \cref{rem:kanmin} a Reedy fibrant simplicial space $K$ is {\it minimal} if $K =\cMin(Y)_0$ (\ref{eq:min ree s}) for some Reedy fibration $Y \to \Delta[0,0,0]$. Let $\uReemin \hookrightarrow \uRee$ be the (Dwyer-Kan equivalent) full subcategory, by \cref{prop:factorization s}.
\end{remone} 

Using \ref{eq:nerve} we obtain a Segal category $N_\Delta\Ree$. We now have the following result.

\begin{theone} \label{the:main simp spaces}
	There is an equivalence between the complete Segal space $\ksS$ and the strict Segal category $N_\Delta\uRee$. Moreover, we have a bijection $\ReeLFib \cong \Hom(-,\ksS)$.
\end{theone}

\begin{proof}
	$\kS$ is a complete Segal space by \cref{prop:spaces s} and the bijection follows from \ref{eq:reelfib}. Hence, we only need to prove the equivalence.
	Let
	\begin{equation} \label{eq:special lfib s}
		\ssL_n = \coprod_{X_0,...,X_n} X_0 \times \Map(X_0,X_1) \times ... \times \Map(X_{n-1},X_n)
	\end{equation}  
	be the trisimplicial set with evident projection map $\pi_1: \ssL \to N_\Delta\Ree$ and notice this map induces a strict pullback $\ssL_n \cong N_\Delta\uRee_n \cong_{\Obj_{\uRee}} \ssL_0$, and so, by applying \cref{rem:strict left} level-wise, the Reedy fibrant replacement over $N_\Delta\uRee$ is a Reedy left fibration. The bijection \ref{eq:reelfib} gives us a functor $N_\Delta\uRee \to \ksS$. We want to prove this is a Dwyer-Kan equivalence. Now, combining \cref{lemma:s vs s min s} and \cref{rem:reemin}, this proof reduces to showing that $\sLift: N_\Delta\uReemin \to \ksS^{min}$ is a Dwyer-Kan equivalence, which follows from establishing the conditions in \cref{prop:segal cat}. 
	
	By \cref{rem:lfib constant s} Reedy left fibrations over $\Delta[0,0]$ are determined by a Reedy fibrant simplicial space and so functor is surjective, meaning we only need to show it is an equivalence of mapping spaces. We can now follow the same steps as in \cref{the:main spaces} to deduce that $\sLift_1: N_\Delta\uReemin \to \ksS^{min}_1$ is a bijection by constructing an explicit inverse.
\end{proof}

\begin{remone} \label{rem:shulman}
	General Reedy fibrant simplicial spaces are too broad to give us $\infty$-categories, however, bisimplicial sets with the Reedy model structure do give us a model topos \cite{rezk2010toposes}, which have been established as models for homotopy type theories \cite{shulman2019inftytoposunivalent}, and so have been studied extensively by Shulman \cite{shulman2015elegantunivalence,shulman2015homotopycanonicity,shulman2017eidiagrams}. In particular bisimplicial sets are the key example of  a model of {\it type theories with shapes} introduced by Riehl and Shulman with the goal of $\infty$-category theory internal to type theories \cite{riehlshulman2017rezktypes}. From this perspective an explicit construction of a complete Segal space for Reedy fibrant simplicial spaces is an important step towards better understanding its properties as a model of type theories. 
\end{remone}

We can now restrict this construction in the following two ways to get the results we wanted. Using the fact that (Segal) coCartesian fibration are local (by \cref{def:reelfib}) we denote $\ksssSet^{\SegcoCart}$ by $\kSeg$ and $\ksssSet^{\coCart}$ by $\kCSS$. We now have the following lemma and theorem, giving us a complete Segal space of (complete) Segal spaces. 

\begin{lemone} \label{lemma:sub css}
 Let $W$ be a complete Segal space. Let $V_{00} \subseteq W_{00}$ be a subset of objects in $W$ closed under equivalences, meaning any object in $W$ equivalent to an element in $V_{00}$ is already in $V_{00}$. Define $V \subseteq W$ as the sub bisimplicial set of $W$ with elements in $\sigma \in V_{nl} \subseteq W_{nl}$ if for all $d:W_{nl} \to W_{00}$, $d(\sigma) \in V_{00}$. Then $V$ is a complete Segal space and $V \to W$ is fully faithful.
\end{lemone}

\begin{proof}
	Let $i: A \to B$ be a trivial cofibration in the complete Segal space model structure and let $f: A \to V$ be a map. Then there exists a map $\hat{f}:B \to W$ lifting $i$. Now, every trivial cofibration in the complete Segal space model structure is surjective on equivalence classes of objects \cite[Lemma 3.54]{rasekh2021yonedad}. Hence, the objects in $W$ that lie in the image of $\hat{f}$ are all equivalent to objects in the image of $f: A \to V \subseteq W$ and hence in $V$ themselves. This proves that the lift factors through $V$ and hence $V$ is a complete Segal space. Finally, for given objects $x,y$, by the definition of $V$ and mapping spaces of complete Segal spaces (\ref{eq:map}), $\Map_V(x,y) \to \Map_W(x,y)$ is the identity and so we are done.
\end{proof}

\begin{theone} \label{the:css}
	We have a diagram of fully faithful functors of complete Segal spaces 
	$$\kCSS \hookrightarrow \kSeg \hookrightarrow \ksS,$$
	with $\kSeg$ having elements Segal spaces and $\kCSS$ complete Segal spaces. Moreover, we have bijections 
	$$\SegcoCart \cong \Hom_{\ssSet}(-,\kSeg),$$
	$$\coCart \cong \Hom_{\ssSet}(-,\kCSS).$$
\end{theone}

\begin{proof}
	The bijections follow directly from \ref{eq:local equiv} and the fact that (Segal) coCartesian fibrations are local (\cref{def:reelfib}). Now, by  \cref{the:main simp spaces}, $\ksS$ is a complete Segal space of bisimplicial sets and, by \cref{lemma:sub css} and \cref{rem:fiberwise cart fib}, $\kSeg \hookrightarrow \ksS$ is a fully faithful inclusion of complete Segal spaces. Here we are using the fact that the Segal condition in \cref{def:reelfib} is by definition up to equivalence and so any Reedy left fibration equivalent to a Segal coCartesian fibration is in fact a Segal coCartesian fibration, proving that the condition in \cref{lemma:sub css} is in fact satisfied. Finally, by \cref{rem:lfib constant s}, the objects in $\kSeg$ are Segal spaces. We can use the same arguments to prove that $\kCSS\to \kSeg$ is a fully faithful functor of complete Segal spaces with $\kCSS$ having objects complete Segal spaces.
\end{proof}

\begin{remone} \label{rem:rfib s fib}
	Combining \cref{rem:rfib fib} and \cref{rem:rfib s} directly implies that the bisimplicial set with $(k,n)$-simplices given by Reedy right fibrations over $\Delta[n,l]$ is precisely the opposite complete Segal space of Reedy fibrant simplicial spaces, $\ksS^{op}$. Moreover, restricting those to (Segal) Cartesian fibrations gives us the opposite complete Segal spaces of (complete) Segal spaces as constructed in  \cref{the:css}, $\kSeg^{op}$ and $\kCSS^{op}$.
\end{remone}

\section{Universal Fibrations} \label{sec:universal}
Up until this point we have constructed various complete Segal spaces that have relevant universal properties in the sense that functors into them correspond to various fibrations over them. This in particular implies the existence of a universal fibration corresponding to the identity map. In this section we want to focus on these universal fibrations. 

We will start with the case for spaces. Denote by $p_{\kssSet}:\kssSet_* \to \kssSet$ the map that corresponds to the identity map under the bijection in \cref{cor:all rep}. Notice, a map $\Delta[n,l] \to \kssSet_*$ corresponds to a map $\sigma:\Delta[n,l] \to \kssSet$ along with a section of the pullback diagram $\sigma^*p_{\kssSet}:\sigma^*\kssSet_* \to \Delta[n,l]$. This means we have a bijection of sets 
\begin{equation} \label{eq:ssset star}
(\kssSet_*)_{nl} \cong \coprod_{p \in \kssSet_{nl}} \Hom_{/\Delta[n,l]}(\id,p),
\end{equation} 
with $p_{\kssSet}$ being the evident projection to $p$. Moreover, the naturality of the bijection  \cref{cor:all rep} implies the following helpful result.

\begin{lemone} \label{lemma:universal pullback}
	The bijection $\All \cong \Hom_{\ssSet}(-,\kssSet)$ is induced by pulling back along $p_{\kssSet}$.
\end{lemone}

Now, for every local class of morphisms $S$, we can obtain $p_{\kssSet^S}:\kssSet_*^S \to \kssSet^S$ by using the bijection \cref{lemma:local bijection}. This map satisfies the following simple, yet useful, lemma that helps us better understand it.

\begin{lemone} \label{lemma:pullback universal fib}
	Let $S$ be a local class of morphisms of bisimplicial sets. Then we have the following pullback square 
	\begin{center}
		\begin{tikzcd}
			\kssSet_*^S \arrow[r, hookrightarrow] \arrow[d] \arrow[dr, phantom, "\ulcorner", very near start] & \kssSet_* \arrow[d] \\
			\kssSet^S \arrow[r, hookrightarrow] & \kssSet
		\end{tikzcd}.
	\end{center}
	In particular, elements in $(\kssSet_*^S)_{nl}$ are morphisms $p:X \to \Delta[n,l]$ that are in $S$ along with a choice of section.
\end{lemone}

 We can now in particular apply this to $S= \LFib$ and deduce that the {\it universal left fibration}, that we denote by $p_{\kS}:\kS_* \to \kS$, has as elements in $(\ksS_*)_{nl}$ diagrams of left fibrations  $\Delta[n,l] \to L$ over $\Delta[n,l]$. We now want to study the representability of this left fibration. Recall that a left fibration $L \to W$ is {\it representable} if there exists an object $x$ in $W$ such that 
 \begin{equation} \label{eq:rep lfib}
  L \simeq W_{x/} = W^{\Delta[1,0]} \times_W \Delta[0,0]. 
 \end{equation} 
 See \cite[Subsection 3.3]{rasekh2017left} for a more detailed analysis of representable left fibrations of bisimplicial sets. 

\begin{theone}\label{the:lfib rep}
	The universal left fibration $p_{\kS}:\kS_* \to \kS$ is a representable left fibration, represented by the terminal object. Moreover the bijection \ref{eq:lfib} is induced by pulling back the universal left fibration $p_{\kS}$.
\end{theone}

\begin{proof}
	The fact that the bijection \ref{eq:lfib} is induced by pulling back $p_{\kS}$ follows directly from \cref{lemma:pullback universal fib} and \cref{lemma:universal pullback}. We now want to prove that $p_{\kS}:\kS_* \to \kS$ is representable and concretely represented by the object $\id_{\Delta[0,0]}$ in $\kS$. By \cite[Theorem 3.55]{rasekh2017left} it suffices to prove that $\kS_*$ has an initial object in the fiber of $p_{\kS}$ over $\id_{\Delta[0,0]}$.

	By \cref{the:main spaces}, we have an equivalence $N_\Delta\uKan \to \kS$, which is induced by a left fibration $\hat{\sL}$ over $\N_\Delta\uKan$, which implies that $\sL \to N_\Delta\uKan$, as constructed in \ref{eq:special lfib}, is equivalent to the homotopy pullback of $p_{\kS}$ along the complete Segal space equivalence $N_\Delta\uKan \to \kS$. Hence, by \cite[Theorem 4.31]{rasekh2017left}, we have a complete Segal space equivalence $\sL\simeq \kS_*$ that takes $(\Delta[0],0)$ to $\id_{\Delta[0,0]}$, where we used the fact that by \ref{eq:special lfib}, $\sL_{00} = \coprod_{X \in \Kan} X_0$, meaning objects in $\sL$ are of the form $(X,x \in X_0)$, where $X$ is a Kan complex.
	
	This implies that in order to finish the proof we only need to observe that $(\Delta[0],\{0\})$ is initial in $\sL$. By definition $\sL_1 = \coprod_{X,Y} X \times \Map_{\uKan}(X,Y)$ and so for an object $(X,x)$, by \ref{eq:map} the mapping space $\Map_{\sL}((\Delta[0],0),(X,x))$ is given via the following pullback
	\begin{center}
		\begin{tikzcd}[row sep=0.4in]
			\Map_{\sL}((\Delta[0],0),(X,x))  \arrow[d, "\cong"] \arrow[r] \arrow[dr, phantom, "\ulcorner", very near start]& \Delta[0] \times \Map_{\uKan}(\Delta[0],X) \arrow[d, "\mathrm{ev}", "\cong"']  \\
			\Delta[0] \arrow[r, "x"] & X
		\end{tikzcd}
	\end{center}
	Now the map on the right hand side is a bijection and so $\Map_{\sL}((\Delta[0],0),(X,x))$ is bijective to $\Delta[0]$ as well, finishing the proof.
\end{proof} 

\begin{remone}
	One of the main results regarding left fibration of bisimplicial sets is that they are always fibrations in the complete Segal space model structure \cite[Corollary 5.11]{rasekh2017left}. Using \cref{the:lfib rep} we can deduce the following result more simply. Indeed, \cref{the:lfib rep} and \ref{eq:rep lfib} imply that $\kS_* \simeq \kS^{\Delta[1,0]} \times_{\kS} \Delta[0,0]$ is a complete Segal space as the complete Segal model structure is Cartesian (see \cite[Lemma 3.43]{rasekh2017left} for a more detailed argument). As a result, $p_{\kS}$ is a Reedy fibration between complete Segal spaces and so a complete Segal fibration \cite[Theorem 7.2]{rezk2001css}. Now, by \cref{the:lfib rep}, every left fibration is a pullback of the complete Segal space fibration $\kS_* \to \kS$ and so a complete Segal space fibration as well. 
\end{remone} 

We now want to generalize \cref{the:lfib rep} to trisimplicial sets. Let $p_{\ksssSet}: \ksssSet_* \to \ksssSet$ be the map that corresponds to the identity map under the bijection \ref{eq:equiv s}. First we want to generalize \ref{eq:ssset star}.

\begin{lemone}
	There is a bijection $$(\ksssSet_*)_{knl} \cong \coprod_{p \in \ksssSet_{nl}} \Hom_{/\Delta[n,l]}(\Delta[k,n,l],p),$$
	meaning an element in $(\ksssSet_*)_{knl}$ corresponds to a choice of trisimplicial set $X \to \Delta[n,l]$ along with a choice of section for the map of bisimplicial sets $X_k \to \Delta[n,l]$.
\end{lemone}

\begin{proof}
	An element in $(\sssSet_*)_{knl}$ is given by a map $\Delta[k,n,l] \to \sssSet_*$, which is given by a map $\Delta[k,n,l] \to \ksssSet$ along with a lift. By \ref{eq:equiv s} an element in $\ksssSet_{knl}$ is given by a map of trisimplicial sets $X \to \Delta[0,n,l]$. Now for such a fixed map, a lift precisely corresponds to a choice of element in $X_{knl}$, which is precisely the data of a section $X_k \to \Delta[n,l]$. 
\end{proof}

This also has the following useful implication.

\begin{lemone} \label{lemma:universal pullback s}
	The bijection $\sAll \cong \Hom_{\ssSet}(-,\ksssSet)$ is induced by pulling back along $p_{\ksssSet}$.
\end{lemone}

Now, for every local class of morphisms $S$, we can obtain $p_{\ksssSet^S}:\ksssSet_*^S \to \ksssSet^S$ by using the bijection \ref{eq:local equiv} and again we have a result analogous to \cref{lemma:pullback universal fib}.

\begin{lemone} \label{lemma:pullback universal fib s}
	Let $S$ be a local class of morphisms of bisimplicial sets. Then we have the following pullback square 
	\begin{center}
		\begin{tikzcd}
			\ksssSet_*^S \arrow[r, hookrightarrow] \arrow[d] \arrow[dr, phantom, "\ulcorner", very near start] & \ksssSet_* \arrow[d] \\
			\ksssSet^S \arrow[r, hookrightarrow] & \ksssSet
		\end{tikzcd}.
	\end{center}
	In particular, elements in $(\ksssSet_*^S)_{knl}$ are morphisms $p:X \to \Delta[n,l]$ that are in $S$ along with a choice of section of $p_k:X_k \to \Delta[n,l]$.
\end{lemone}

We can now apply this lemma to the case $S=\ReeLFib$ to obtain the universal Reedy left fibration $p_{\ReeLFib}:\ksS_* \to \ksS$. We want to observe that this Reedy left fibration is representable, as defined in \cite[Definition 2.1]{rasekh2017cartesian}.

\begin{theone}\label{the:ree lfib rep}
	The universal Reedy left fibration $p_{\ksS}:\ksS_* \to \ksS$ is a representable Reedy left fibration represented by the cosimplicial object $\Delta[\bullet,0]$. Moreover the bijection \ref{eq:reelfib} is induced by pulling back the universal left fibration $p_{\kS}$.
\end{theone}

\begin{proof}
	The fact that the bijection is induced by pullback follows from combining \cref{lemma:pullback universal fib s} and \cref{lemma:universal pullback s}. Now, in order to prove that $p_{\ksS}: \ksS_*\to \ksS$ is representable, by \cite[Lemma 4.1]{rasekh2017cartesian}, it suffices to prove that the left fibration $(p_{\ksS})_k: (\ksS_*)_k\to \ksS$ is represented by $\Delta[k,0]$ for all $k \geq 0$. From here one we can follow the steps of the proof of \cref{the:lfib rep}.
	
	We need to show that $(\ksS_*)_k$ has an initial object over $\Delta[k,0]$. By \cref{the:main simp spaces}, we have an equivalence $N_\Delta\uRee \to \ksS$, which is induced by the Reedy left fibration $\hat{\ssL}$ over $N_\Delta\uRee$, which implies that the Reedy left fibration $\ssL \to N_\Delta\uRee$, as constructed in \ref{eq:special lfib s}, is the homotopy pullback of the Reedy left fibration $\ksS_* \to \ksS$, and, in particular, the left fibration $\hat{\ssL}_k \to N_\Delta\uRee$ is the homotopy pullback of the left fibration $(\ksS_*)_k \to \ksS$. Hence, by \cite[Theorem 4.31]{rasekh2017left}, we have a complete Segal space equivalence $\ssL_k \simeq (\ksS_*)_k$ that takes $(\Delta[k,0],\id_{k})$ to $\id_{\Delta[k,0]}$. Here we used the fact that by \ref{eq:special lfib s}, $\ssL_{k00} = \coprod_{X \in \Ree} X_{k0}$, meaning objects in $\ssL_k$ are of the form $(X,x \in X_{k0})$, where $X$ is a Reedy fibrant simplicial space.
	
	This implies that in order to finish the proof we only need to observe that $(\Delta[k,0],\id_{k})$ is initial in $\ssL_k$. By definition $\ssL_{k1} = \coprod_{X,Y \in \Ree} X_k \times \Map_{\uRee}(X,Y)$ and so for an object $(X,x \in X_{k0})$, by \ref{eq:map} the mapping space $\Map_{\ssL_k}((\Delta[k,0],\id_k),(X,x))$ is given via the following pullback
	\begin{center}
		\begin{tikzcd}
			\Map_{\ssL_k}((\Delta[k,0],\id_k),(X,x))  \arrow[d, "\cong"] \arrow[r] \arrow[dr, phantom, "\ulcorner", very near start]& \Delta[k,0]_k \times \Map_{\uRee}(\Delta[k,0],X) \arrow[d, "\mathrm{ev}", "\cong"']  \\
			\Delta[0] \arrow[r, "(\id_k \comma x)"] & \Delta[k,0]_k \times X_k
		\end{tikzcd}
	\end{center}
	Now, the map on the right hand side is a bijection, as $\Map_{\uRee}(\Delta[k,0],X) \cong X_k$ by \ref{eq:map reedy}. As a result, the mapping space  $\Map_{\ssL_k}((\Delta[k,0],\id_k),(X,x))$ is bijective to $\Delta[0]$ as well, finishing the proof.
\end{proof} 

We now use \cref{the:css} to get the universal Segal coCartesian fibration, that we denote by $p_{\kSeg}:\kSeg_* \to \kSeg$ and the universal coCartesian fibration, that we denote by $p_{\kCSS}:\kCSS_* \to \kCSS$. \cref{lemma:pullback universal fib s} and \cref{the:ree lfib rep} now immediately give the following result.

\begin{corone} \label{cor:universal fib css}
	The universal (Segal) coCartesian fibration $p_{\kCSS}$ ($p_{\kSeg}$) is represented by $\Delta[\bullet,0]$. Moreover, we have pullback squares
	\begin{center}
		\begin{tikzcd}
			\kCSS_* \arrow[d, "p_{\kCSS}"] \arrow[r, hookrightarrow] \arrow[dr, phantom, "\ulcorner", very near start] & \kSeg_* \arrow[r, hookrightarrow] \arrow[d, "p_{\kSeg}"] \arrow[dr, phantom, "\ulcorner", very near start] & \ksS_* \arrow[d, "p_{\ksS}"] \\
			\kCSS \arrow[r, hookrightarrow] & \kSeg \arrow[r, hookrightarrow] & \ksS
		\end{tikzcd}.
	\end{center}
	Finally, pulling back along the universal fibrations induces bijections
	$$\Hom_{\ssSet}(-,\kSeg) \cong \SegcoCart,$$
	$$\Hom_{\ssSet}(-,\kCSS) \cong \coCart.$$
\end{corone}

\begin{remone} \label{rem:rep cart fib}
	The representability of the universal left fibration is well-established (and has for example been studied in \cite[Subsection 5.2]{cisinski2019highercategories},) however, the representability of the the universal coCartesian fibration is a more modern phenomena and can be found in \cite[Subsection 4.2]{rasekh2017cartesian} and \cite[Example 3.26]{stenzel2020comprehension}. The representability of the universal Reedy left fibration was not studied before.
\end{remone}

\section{Comparison with Quasi-Categories}
Up until here we constructed various complete Segal spaces of spaces and (complete) Segal spaces. In this last section we want to use the fact that we can translate between complete Segal spaces and quasi-categories, another important model of $(\infty,1)$-categories, to construct quasi-categories of spaces and (complete Segal) spaces. This requires us to review left fibrations of simplicial sets \cite{joyal2008theory,lurie2009htt,cisinski2019highercategories}, the translation results between quasi-categories and complete Segal spaces due to Joyal and Tierney \cite{joyaltierney2007qcatvssegal} and their generalization to left fibrations in \cite[Appendix B]{rasekh2017left}. 

\begin{defone} \label{def:lfib sset}
	A {\it left fibration} of simplicial sets is map that satisfies the right lifting property with respect to horn inclusions $\Lambda[n]_i \hookrightarrow \Delta[n]$, for $0\leq i < n$.
\end{defone}
Left fibrations of simplicial sets can be translated to left fibrations of bisimplicial sets (\cref{def:left fib}) and vice versa. Let $i_1^*:\ssSet \to \sSet$ be the functor that takes a bisimplicial set $X_{\bullet\bullet}$ to the simplicial set $X_{\bullet 0}$ \cite[Section 4]{joyaltierney2007qcatvssegal}. Moreover, let $t^!:\sSet \to \ssSet$ be the functor that takes a simplicial set $X$ to the bisimplicial set $t^!X_{nl} = \Hom_{\sSet}(\Delta[n] \times N(I[l]),X)$ \cite[Theorem 2.12]{joyaltierney2007qcatvssegal}. Here $I[l]$ is the groupoid with $l+1$ objects and a unique morphism between any two objects. 

These two functors are both right adjoints of Quillen equivalences \cite[Theorem 4.11, Theorem 4.12]{joyaltierney2007qcatvssegal}, which in particular has the following implications:

\begin{lemone} \label{lemma:qcat vs seg jt}
	$i_1^*$ takes complete Segal spaces to quasi-categories and $t^!$ takes quasi-categories to complete Segal spaces. Moreover, $i_1^*t^!:\sSet \to \sSet$ is the identity map and $t^!i_1^*: \ssSet \to \ssSet$ is equivalent to the identity. Finally, $i_1^*$ reflects equivalences between complete Segal spaces.
\end{lemone}

 These results have been generalized in \cite[Appendix B]{rasekh2017left} to a comparison between left fibrations of simplicial sets and bisimplicial sets, giving us the following valuable result.
 
 \begin{lemone} \label{lemma:qcat vs seg left}
 	$i_1^*$ takes left fibrations of bisimplicial sets to left fibrations of simplicial sets and $t^!$ takes left fibrations of simplicial sets to left fibrations of bisimplicial sets.
 \end{lemone}

We now use the ability to translate between quasi-categories and complete Segal spaces to construct additional $(\infty,1)$-categories of spaces. First of all we can apply $i_1^*$ to the complete Segal space $\kS$ (\cref{the:main spaces}) to obtain the following result.

\begin{corone} \label{cor:ione qcat}
	$i_1^*\kS$ is a quasi-category of spaces with $i_1^*\kS_n$ given by left fibrations of bisimplicial sets over $\Delta[n,0]$.
\end{corone}

We now want to illustrate how we can use left fibrations of simplicial sets internally to construct a quasi-category of spaces, using analogous steps to \cref{sec:spaces}. Let $\kS_{\QCat}$ be the simplicial set with $(\kS_{\QCat})_n$ given by left fibrations of simplicial sets over $\Delta[n]$ (where we are using the translation to functors as given in \cref{lemma:bijection} to take care of functoriality). 

Now, by \cref{lemma:qcat vs seg left}, $t^!$ preserves left fibrations and moreover we have $t^!(\Delta[n]) = \Delta[n,0]$. Hence $t^!$ induces a morphism of quasi-categories 
\begin{equation} \label{eq:kt}
	\kT:\kS_{\QCat} \to i_1^*\kS,
\end{equation}
that takes a left fibration of simplicial sets $L \to \Delta[n]$ to $t^!L \to t^!\Delta[n] = \Delta[n,0]$. Similarly, by \cref{lemma:qcat vs seg left}, $i_1^*$ also preserves left fibrations and $i_1^*(\Delta[n,0]) = \Delta[n]$ and so $i_1^*$ similarly induces a morphism of quasi-categories
\begin{equation} \label{eq:ki}
 \kI: i_1^*\kS \to \kS_{\QCat},	
\end{equation}
that takes a left fibration of bisimplicial sets $L \to \Delta[n,0]$ to $i_1^*(L) \to i_1^*(\Delta[n,0]) = \Delta[n]$. We now have the following result.

\begin{theone} \label{the:qcat v segal}
	The maps $\kT:\kS_{\QCat} \to i_1^*\kS$ and $\kI: i_1^*\kS \to \kS_{\QCat}$ are inverses of quasi-categories.
\end{theone} 

\begin{proof}
	First we prove $\kS_{\QCat}$ is a quasi-category. By \cref{lemma:qcat vs seg jt}, $t^!\circ i_1^*$ is the identity and so $\kI\kT$ is the identity as well, meaning $\kS_{\QCat}$ is a retract of the quasi-category $i_1^*\kS$ (\cref{cor:ione qcat}) and so a quasi-category as well. We now move to prove $\kT$ and $\kI$ are inverses of quasi-categories. As $\kI\kT$ is the identity we only need to show $\kT\kI:i_1^*\kS \to i_1^*\kS_{\CSS}$ is equivalent to the identity. By \cref{lemma:qcat vs seg jt} and \cref{lemma:s vs s min}, the two morphisms $\i_1^*\kS^{min} \to i_1^*\kS \to \i_1^*\kS^{min}$ are equivalences of quasi-categories, hence it suffices to prove that $\kT\circ \kI \circ (i_1^*\kMin):i_1^*\kS^{min} \to i_1^*\kS^{min}$ is equal to the identity. 
	
	Let $L \to \Delta[n,0]$ be a left fibration. By \cref{lemma:qcat vs seg jt}, there is an equivalence of complete Segal spaces $t^!i_1^*L \to L$ over $\Delta[n,0]$, which implies they are equivalent left fibrations \cite[Theorem 5.11]{rasekh2017left} and so $\cMin(t^!i_1^*L)$ and $\cMin(L)$ are equal (\cref{lemma:reedy equiv strict}) finishing the proof. 
\end{proof}

\begin{remone} \label{rem:cisinski}
	The elements in the quasi-category $\kS_{\QCat}$ are precisely left fibrations over $\Delta[n]$. Hence this constructions coincides with the construction of the quasi-category of spaces by Cisinski \cite[Theorem 5.2.10, Corollary 5.4.7]{cisinski2019highercategories} and hence gives us an independent proof thereof.
\end{remone}

We can take the opposite route to \cref{cor:ione qcat} to get the following result.

\begin{corone} \label{cor:tshriek css}
	$t^!\kS_{\QCat}$ is a complete Segal space of spaces with $(t^!\kS_{\QCat})_{nl}$ given by left fibrations of simplicial sets over $\Delta[n] \times N(I[l])$.
\end{corone}

We can now use this to generalize to quasi-categories of complete Segal spaces. 

\begin{corone}
	$i_1^*\ksS$ is a quasi-category with $i_1^*\ksS_n$ given by Reedy left fibrations of trisimplicial sets over $\Delta[0,n,0]$. Moreover, we have inclusions of quasi-categories 
	$$i_1^*\kCSS \hookrightarrow i_1^*\kSeg \hookrightarrow i_1^*\ksS,$$
	where $i_1^*\kCSS$ and $i_1^*\kSeg$ are characterized analogously.
\end{corone}

There is also a result analogous \cref{the:qcat v segal} by changing the elements of the complete Segal space of complete Segal spaces. A Reedy fibration of bisimplicial sets $L \to \Delta[0,n]$ is a Reedy left fibration if for all $k\geq 0$ the restricted map $L_k \to \Delta[0,n]_k = \Delta[n]$ (\cref{not:delta dash dash}) is a left fibration of simplicial sets. Let $st^!:\ssSet \to \sssSet$ be defined as $\Fun(\DD^{op},t^!)$ and similarly, let $si_1^*= \Fun(\DD^{op},i_1^*):\ssSet \to \sssSet$. We now have the following fact about $si_1^*$ and $st^!$.

\begin{lemone} \label{lemma:qcat vs seg jt s}
	$si_1^*$ takes Reedy left fibrations of trisimplicial sets to Reedy left fibrations of bisimplicial sets and $st^!$ takes Reedy left fibrations of bisimplicial sets to Reedy left fibrations of trisimplicial sets. Moreover, $si_1^*st^!$ is the identity and $st^!si_1^*$ takes a Reedy left fibration of trisimplicial sets to an equivalent one.
\end{lemone}
See \cite[Theorem 1.35]{rasekh2021cartfibmarkedvscso} for further details about these functors. Using this result we can generalize the maps of quasi-categories $\kT,\kI$ to maps
\begin{equation}\label{eq:skt}
	\begin{aligned} 
		\skT:\ksS_{\QCat} \to i_1^*\ksS,  \\
		\skI: i_1^*\ksS \to \ksS_{\QCat},
	\end{aligned}
\end{equation} 
and following the same steps of the proof of \cref{the:qcat v segal} (this time with the construction of minimal Reedy left fibrations given in \ref{eq:min ree s}) gives us the following result.

\begin{theone} \label{the:qcat v segal s}
	The maps $\skT:\ksS_{\QCat} \to i_1^*\ksS$ and $\skI: i_1^*\ksS \to \ksS_{\QCat}$ are inverses of quasi-categories.
\end{theone}

Finally, we can use the same conditions as in \cref{def:reelfib} to define (Segal) coCartesian fibrations of bisimplicial sets over $\Delta[0,n]$. 
As a result we get the subsimplicial sets $\kCSS_{\QCat} \hookrightarrow \kSeg_{\QCat} \hookrightarrow \ksS_{\QCat}$ of (Segal) coCartesian fibrations. Now, by \cite[Theorem 1.35]{rasekh2021cartfibmarkedvscso}, $st^!$ and $si_1^*$ preserve (and in fact reflect) (Segal) coCartesian fibrations, meaning, our functors $\skT$ and $\skI$ restrict appropriately, giving us the following result.

\begin{corone} \label{cor:big diagram}
	We have the following diagram of quasi-categories, where the horizontal maps are inclusions and the vertical maps equivalences.
	\begin{center}
		\begin{tikzcd}
			i_1^*\kCSS \arrow[r, hookrightarrow] \arrow[d, "\skI", "\simeq"'] & i_1^*\kSeg \arrow[r, hookrightarrow] \arrow[d, "\skI", "\simeq"'] & i_1^*\ksS \arrow[d, "\skI", "\simeq"'] \\
			\kCSS_{\QCat} \arrow[r, hookrightarrow] \arrow[d, "\skT", "\simeq"'] & \kSeg_{\QCat} \arrow[r, hookrightarrow]  \arrow[d, "\skT", "\simeq"'] & \ksS_{\QCat}  \arrow[d, "\skT", "\simeq"'] \\
			i_1^*\kCSS \arrow[r, hookrightarrow] & i_1^*\kSeg \arrow[r, hookrightarrow] & i_1^*\ksS \\
		\end{tikzcd}
	\end{center}
\end{corone}

Finally, we have the following result analogous to \cref{cor:tshriek css}

\begin{corone} \label{cor:tshriek css s}
	$t^!\kCSS_{\QCat}$ is a complete Segal space of complete Segal spaces spaces with $(t^!\kCSS_{\QCat})_{nl}$ given by coCartesian fibrations of bisimplicial sets over $\Delta[0,n] \times N(I[l])$.
\end{corone}

\bibliographystyle{alpha}
\bibliography{main}

\end{document}